\documentclass[amstex,12pt]{amsart}
\usepackage{geometry}
\usepackage{amsmath, mathtools,amssymb,amsthm,mathrsfs,color,lineno,paralist,graphicx,float}
\usepackage[linktocpage,backref=page,colorlinks,
linkcolor=green,
anchorcolor=red,
citecolor=red,
]{hyperref}
\usepackage{booktabs,tabularx}
\usepackage{color}
\usepackage[T1]{fontenc}
\usepackage[utf8]{inputenc}
\usepackage{mathrsfs}
\usepackage{cleveref}
\usepackage{mathscinet}
\usepackage[normalem]{ulem}
\usepackage{calc}
\newtheorem{thm}{Theorem}[section]
\theoremstyle{plain}
\theoremstyle{definition}
\newtheorem{proposition}[thm]{Proposition}

\newtheorem{cor}[thm]{Corollary}
\newtheorem{lemma}[thm]{Lemma}
\newtheorem{definition}[thm]{Definition}
\newtheorem{claim}[thm]{Claim}
\newtheorem{remark}[thm]{Remark}

\numberwithin{equation}{section}

\AtBeginEnvironment{thm}{\crefalias{thm}{theorem}}
\AtBeginEnvironment{lemma}{\crefalias{thm}{lemma}}
\AtBeginEnvironment{cor}{\crefalias{thm}{corollary}}
\AtBeginEnvironment{proposition}{\crefalias{thm}{proposition}}
\AtBeginEnvironment{claim}{\crefalias{thm}{claim}}

\DeclareMathOperator{\dist}{dist}

\newcommand{\R}{\mathbb R}

\newcommand{\Id}{\operatorname{Id}}

\definecolor{lgray}{gray}{0.9}

\def\leq{\leqslant}
\def\geq{\geqslant}

\def\hat{\widehat}

\allowdisplaybreaks
\renewcommand*{\backref}[1]{}
\renewcommand*{\backrefalt}[4]{\quad \tiny
  \ifcase #1 (\textbf{NOT CITED.})
  \or    (Cited on page~#2.)
  \else   (Cited on page~#2.)
  \fi}
\makeatletter

\title[Uniqueness of SRB measures]{Uniqueness of SRB measures for stochastic differential equations}

\date{\today}

\author[Z. Lian]{Zeng Lian}
\address[Z. Lian]
{	School of Mathematics,
Sichuan University,
Chengdu, Sichuan, 610016, China} \email{lianzeng@scu.edu.cn, zenglian@gmail.com}

\author[J. Zhang]{Jianhua Zhang}
\address[J. Zhang]{School of Mathematical Sciences, 
Chengdu University of Technology, Chengdu, Sichuan, 610059, China}
\email{leapforg@mail.ustc.edu.cn}

\begin{document}

\begin{abstract}
For a dissipation stochastic differential equation (abbreviated as SDE),  exponential  ergodicity   has been investigated a lot in a wide range of settings.  In this paper, we  prove that  the  Markovian measure $\mu$ determined by the unique  stationary measure  is  also the unique invariant measure which satisfies the entropy formula. Consequently,  the existence of an SRB measure is completely determined by the sign of the sample top Lyapunov exponent of $\mu$: if $\chi_\mu^1>0$, then $\mu$ is the  unique SRB measure, whereas if $\chi_\mu^1\leq0$, then the SDE admits no  SRB measure. As an application, we establish a  dichotomy of SRB measures for the 2D additive-noise stochastic Hopf normal form within a shear parameter domain.  
  \end{abstract}
\maketitle



\section{Introduction}

SRB measures,  named after Sinai, Ruelle and Bowen, were  discovered   in the context of uniformly hyperbolic attractors \cite{MR399421,MR415683,MR380889}.
Let $\varphi$ be a $C^2$ diffeomorphism  on a closed Riemannian manifold.  In general, an invariant   measure $\mu$ of $\varphi$ is called   an SRB measure if (i) $(\varphi,\mu)$ has  positive Lyapunov exponents; (ii) the conditional measures of $\mu$ on unstable manifolds are absolutely continuous with respect to the volume  measures on these manifolds. A fundamental connection between Lyapunov exponents and metric entropy is provided by  the entropy formula.  Pesin \cite{MR466791} proved that, if an invariant  measure $\mu$ is absolutely continuous
with respect to Riemannian volume measure, then
\begin{equation}
\label{eq:introduction-pesin-formula}
 h_\mu(\varphi)
 =\int\sum_{i=1}^d\chi^i_{\mu}(x)^+\,d\mu(x),
\end{equation}
where  the Lyapunov exponents $\chi^i_{\mu}$ are counted with multiplicities and $a^+=\max\{a,0\}$.  Ledrappier and Young in \cite{MR819556} subsequently characterized SRB
measures, in the presence of positive Lyapunov exponents, in terms of
the entropy formula \eqref{eq:introduction-pesin-formula}.  This relationship between entropy, Lyapunov
exponents, and dimension  has since
been extended to a variety of  settings (\cite{MR968818,MR1646606,MR1676734,MR4794090,MR3595937,MR3038055}).

The theory of SRB measures  has also been developed for random dynamical systems
\cite{MR3925383,MR968818,MR1646606,MR1369243}.  In the present paper, we study the random dynamical system  generated by the  following SDE on $\mathbb{R}^d$
\begin{equation}
\label{26-6-8-1611}
 dX_t=b(X_t)dt+\sum_{i=1}^m\sigma_i(X_t) dW_t^i.
\end{equation}
Under suitable assumptions,  equation \eqref{26-6-8-1611} generates a smooth  random dynamical
system. There is a substantial  of literatures to investigate the exponential ergodicity  (for example, \cite{MR2259251,MR1846802,MeynTweedie2009, MR2857021}).  In particular,  there exists  a Lyapunov-Foster function $V:\mathbb{R}^d\to [1,\infty)$, a unique stationary measure $\gamma$ on $\mathbb{R}^d$ and a constant $\lambda>0$ such that the Markov semigroup $\mathcal{Q}_t$ satisfies 
\begin{align}
\|\mathcal{Q}_tf-\gamma(f)\|_V\lesssim e^{-\lambda t}\|f\|_{V}
\end{align}
for any  Borel function $f$ with $\|f\|_V:=\sup_{x\in\mathbb{R}^d}\frac{|f(x)|}{V(x)}<\infty$. This  stationary measure determines  a  unique Markovian measure which is an invariant measure   for the  random dynamical system  (see \Cref{26-6-9-1647}).  It is clear that  the exponential ergodicity implies the unique  stationary measure has a density\footnote{Let $\varphi_t$ be the stochastic flow  generated by the equation \eqref{26-6-8-1611}.  Choose a compactly supported probability measure
$\eta\ll\mathrm{Leb}$.  For any zero-Lebesgue measure Borel set $A$,  we have 
$$\int_{\mathbb{R}^d}\mathcal Q_t\mathbf{1}_A\,d\eta=\int_{\Omega}\int_{\mathbb{R}^d}\mathbf{1}_A\circ\varphi_t\,d\eta\,d\mathbb{P}=0.$$ The  exponential ergodicity  therefore yields
\[
\gamma(A)=\left|\int_{\mathbb{R}^d}\bigl(\gamma(A)-\mathcal Q_t\mathbf1_A(x)\bigr)\,d\eta(x)\right|
\lesssim  e^{-\lambda t}\eta(V)\rightarrow0.
\]
Hence $\gamma\ll\mathrm{Leb}$.}.  Consequently, the  Markovian measure satisfies the  entropy formula \cite{MR968818,MR1369243,MR1646606,MR3175262}. Recently, positivity of the top sample  Lyapunov exponent\footnote{Here the top Lyapunov exponents is taken with respect to the unique Markovian measure.} has
been proved for several concrete stochastic models \cite{ChemnitzEngel2023,MR4372219,MR4733339,breden2024rigorous}.  Together with the characterization of entropy formula,  this unique Markovian measure  is an SRB measure.  

It is natural to ask  whether  equation  \eqref{26-6-8-1611} admits another SRB measures.  In this paper, we discover  a  rigidity phenomenon:  the unique Markovian measure is also the  unique invariant measures which satisfies the entropy formula.  To the authors' knowledge, it is a completely new  rigidity phenomenon. The proof combines two ingredients.  First, we develop the multiplicative mean ergodic theorem for a class of Feynman–Kac semigroups.  Second, we establish a weak version Katok's  formula for non-ergodic measure in non-compact phase spaces and a uniform volume estimate for Bowen balls on the  Pesin blocks in non-compact phase space. As an application, we obtain a   dichotomy  of SRB measures  for the 2D additive-noise stochastic Hopf normal form for the shear parameter.

To ensure that  \eqref{26-6-8-1611} generates a well-defined random dynamical system, we impose the following hypotheses: 
\begin{enumerate}[(H1)]
\item\label{26-6-10-1501} $b\in C^{\infty}(\mathbb{R}^d)$ and $\sigma_i\in C_b^{\infty}(\mathbb{R}^d)$ for $1\leq i\leq m$;
\item\label{26-6-10-1502}  there exist  $\alpha>0,\beta> 0$, $C>1$ and  $p>2$ such that  for any $x\in\mathbb{R}^d$
$$\langle x, b(x) \rangle\leq -\alpha |x|^{p}+\beta,\quad\text{and}\quad |b(x)|+|\nabla b(x)|+|\nabla^2 b(x)|\leq C(|x|^p+1) ;$$
\item\label{26-6-10-1503}  the  diffusion matrix $a(x):=\sum_{i=1}^m\sigma_i(x)\sigma_i(x)^T$  is uniformly elliptic.
\end{enumerate}

Let
\[
 \Omega:=C_0(\R;\R^m)
 :=\{\omega\in C(\R;\R^m):\omega(0)=0\}
\]
be endowed with the compact-open topology,  $\mathscr F$ be its Borel $\sigma$-algebra, and  $\mathbb P$ be the two-sided Wiener measure on its Borel $\sigma$-algebra.  The coordinate process $(W^1_t,\dots, W_t^m)_{t\in\mathbb{R}}$ is an $m$-dimensional Brownian motion.  Now, we define the following two-parameter family of $\sigma$-algebras by 
$$\mathscr{F}_s^t:=\sigma\{ W_{u}-W_v:  s\leq v\leq u\leq t\}$$
for $s\leq t$. For any $t\in\mathbb{R}$,  define the shift transformation  $\theta_t: \Omega\to\Omega$ by 
$$\theta_t(\omega)=\omega(\cdot +t)-\omega(t).$$
Therefore, $(\Omega,\mathscr{F},\mathbb{P},(\theta_t)_{t\in\mathbb{R}})$ is an ergodic metric dynamical system (for example, see \cite[Appendix A]{A}).  By \cite[Chapter~IV]{MR1472487}, together with  \cite[Theorem~5.1]{MR1305784} and the Lyapunov estimates following from \hyperref[26-6-10-1501]{(H1)}--\hyperref[26-6-10-1503]{(H3)},  equation \eqref{26-6-8-1611} generates a $C^2$ \emph{random dynamical system} (abbreviated as RDS).  In particular, there is a Borel measurable mapping 
$\varphi:\mathbb{R}_+\times\Omega\times\mathbb{R}^d\to\mathbb{R}^d$
which   satisfies the following properties
\begin{itemize}
\item  for any $x\in\mathbb{R}^d$, $\{\varphi(t,\cdot, x)\}_{t\geq0}$ is the strong solution  of  equation \eqref{26-6-8-1611}.
\item  for $\mathbb{P}$-a.s. $\omega\in\Omega$  we have  $\varphi(t+s,\omega,x)=\varphi(t,\theta_s\omega, \varphi(s, \omega,x))$ whenever $s,t\in\mathbb{R}_+$ and $x\in\mathbb{R}^d$;
\item for $\mathbb{P}$-a.s. $\omega\in\Omega$, the mapping $\varphi(t, \omega):=\varphi(t,\omega,\cdot):\mathbb{R}^d\to\mathbb{R}^d$ is a $C^2$ injective local diffeomorphism. 
\end{itemize}
The associated skew-product semiflow  $\Phi:=(\Phi_t)_{t\in\mathbb{R}_+}$ is defined  as 
\begin{align}
\label{26-6-9-1558}
\Phi_t: \Omega\times\mathbb{R}^d\to\Omega\times\mathbb{R}^d,\quad (\omega,x)\mapsto (\theta_t\omega, \varphi(t,\omega,x)).
\end{align}
In the remainder of this paper,  unless otherwise stated, the notation  $\varphi$ refers  to the RDS generated by  equation \eqref{26-6-8-1611}.

As shown in the proof of \Cref{25-6-12-1553}, equation \eqref{26-6-8-1611} admits a compact random set
that pullback absorbs every bounded non-random subset of $\mathbb R^d$. By \cite[Theorem~3.11]{MR1305587}, it therefore admits a compact \emph{pullback random attractor}. That is,  there exists a set-valued mapping $\mathfrak{A}: \Omega\to 2^{\mathbb{R}^d}$ satisfying
\begin{enumerate}[(i)]
\item  for any $x\in\mathbb{R}^d$,  $\omega\mapsto \text{dist}(x,\mathfrak{A}(\omega))$ is $\mathscr{F}$-measurable;
\item  for $\mathbb{P}$-a.s. $\omega\in\Omega$,  $\varphi(t, \omega, \mathfrak{A}(\omega))=\mathfrak{A}(\theta_t\omega)$ whenever $t\in\mathbb{R}_+$, and $\mathfrak{A}(\omega)$ is compact;
\item  for any bounded subset $S$ of $\mathbb{R}^d$ one has that for $\mathbb{P}$-a.s. $\omega\in\Omega$
$$\lim_{t\to\infty}\text{dist}(\varphi(t,\theta_{-t}\omega, S), \mathfrak{A}(\omega))=0,$$
 where $\text{dist}(\cdot,\cdot)$ is the Hausdorff semi-distance\footnote{It means that $\text{dist}(A,B)=\sup_{x\in A}\inf_{y\in B}|x-y|$ for the subsets  $A,B$ of $\mathbb{R}^d$.}. 
\end{enumerate}

\subsection{Invariant measures}
In this subsection, we mainly review the concept of invariant measures of RDS and  state the relationship between stationary measures and invariant measures.

The   RDS $\varphi$ on $\mathbb{R}^d$  driven  by    $(\Omega,\mathscr{F},\mathbb{P}, (\theta_t)_{t\in\mathbb{R}})$ generates a  Markovian semigroup
$(\mathcal{Q}_t)_{t\geq 0}$ defined  as 
$$\mathcal{Q}_t(f)=\int_{\Omega} f\circ\varphi(t,\omega,\cdot)\,d\mathbb{P}(\omega),$$
where $f:\mathbb{R}^d\to\mathbb{R}$ is a bounded measurable function.  A Borel probability measure $\gamma$ on $\mathbb{R}^d$ is called a \emph{stationary measure} of Markov semigroup $(\mathcal{Q}_t)_{t\geq 0}$ if 
for any bounded measurable function $f:\mathbb{R}^d\to\mathbb{R}$ one has that 
$$\int_{\mathbb{R}^d}f\,d\gamma= \int_{\mathbb{R}^d}\mathcal{Q}_t(f)\,d\gamma.$$
 A stationary measure $\gamma$  is  called  an \emph{ergodic stationary measure}  if  any $\mathcal{Q}$-invariant  bounded measurable function  is a  constant up to  a $\gamma$-null measure set.  

Now, we recall the definition and properties of invariant measures of RDS. 
\begin{definition}
 A  Borel probability measure $\nu$ on $\Omega\times\mathbb{R}^d$ is called  a  $\varphi$-\emph{(ergodic) invariant measure}  if it satisfies the following properties:
    \begin{enumerate}[(a)]
        \item $\nu$ is invariant (resp. ergodic)  with respect to  the skew-product measurable dynamical system  $\Phi$ defined in \eqref{26-6-9-1558}; 
        \item $(\pi_{\Omega})_*\nu = \mathbb{P}$, where $\pi_{\Omega}$ denotes the projection from $\Omega \times\mathbb{R}^d$ to $\Omega$. 
    \end{enumerate}
    Such a pair $(\varphi, \nu)$ is called a (resp. ergodic) \emph{measure-preserving RDS}.
\end{definition}
For the sake of convenience,  we denote by $\mathcal{I}(\varphi)$  the collection of all invariant measures of  RDS $\varphi$.  In fact, the  Borel  probability measure $\nu$ on $\Omega\times \mathbb{R}^d$ with $(\pi_{\Omega})_*\nu=\mathbb{P}$  is characterized by its sample measures $\{\nu_{\omega}\}_{\omega\in\Omega}$ on $\mathbb{R}^d$  with the following properties 
  \begin{itemize}
    \item for any bounded measurable $f: \mathbb{R}^d\to\mathbb{R}$,  the mapping $\omega\mapsto \int_{\mathbb{R}^d}f\,d\nu_{\omega}$ is $\mathscr{F}$-measurable;
    \item for any bounded measurable $g:\Omega\times \mathbb{R}^d\to\mathbb{R}$, one has that 
    $$\int_{\Omega\times \mathbb{R}^d}g\,d\nu=\int_{\Omega}\int_{\mathbb{R}^d}g\,d\nu_{\omega}d\mathbb{P}(\omega).$$
    \end{itemize}
    It  is written symbolically $\nu=\int_{\Omega}\nu_{\omega}\,\mathrm{d}\mathbb{P}(\omega)$, usually called as the disintegration of $\nu$. The reader can see the details in  \cite[Proposition 3.6]{MR1993844}. 
\begin{lemma}[{\cite[Theorem 1.4.5]{A}}]    
A Borel probability measure $\nu$  on $\Omega\times\mathbb{R}^d$ is an invariant measure of the  RDS $\varphi$ if and only if $(\pi_{\Omega})_*\nu=\mathbb{P}$ and   for all $t\geq 0$
$$\big(\varphi(t,\omega)\big)_*\nu_{\omega}=\nu_{\theta_t\omega}$$ 
for $\mathbb{P}$-a.s. $\omega\in\Omega$.
\end{lemma}

    An invariant measure $\nu=\int_{\Omega}\nu_{\omega}\,d\mathbb{P}(\omega)$ is called  a \emph{Markovian measure} if  $\omega\mapsto\nu_{\omega}$ is $\mathscr{F}^0_{-\infty}$-measurable.  Finally, we state the relationship between Markovian  measures and stationary measures. The reader can see the details in  \cite[Theorem 1.7.2]{A}.
\begin{proposition}
\label{26-6-9-1647}
Let $\gamma$ be a stationary  measure of the RDS $\varphi$ on $\mathbb{R}^d$.  Then  there is a unique Markovian   measure $\mu$  defined by 
    $$\omega\mapsto\mu_{\omega}=\lim_{t\to\infty}(\varphi(t, \theta_{-t}\omega))_*\gamma.$$
 Conversely, $\gamma=\mathbb{E}[\mu_{\cdot}]$.   Furthermore, $\mu$ is ergodic if and only if $\gamma$ is  ergodic. 
\end{proposition}

Under \hyperref[26-6-10-1501]{(H1)}-\hyperref[26-6-10-1503]{(H3)}, we  can apply  the weighted Harris theorem  to establish the exponential ergodicity of equation \eqref{26-6-8-1611} (see  \Cref{clm:weighted-mixing}).   From now on,  unless otherwise stated, we always use  $\gamma$ and $\mu$ to denote   the unique stationary measure and Markovian measure of  equation \eqref{26-6-8-1611}, respectively. 
\subsection{Metric entropy and SRB measures}
In this subsection, we  recall   the definitions of  metric  entropy  and SRB measures for RDS.

Let $(\varphi,\nu)$ be a measure-preserving RDS on $\mathbb{R}^d$ driven by  $(\Omega,\mathscr{F},\mathbb{P}, (\theta_t)_{t\in\mathbb{R}})$.  For any finite  Borel measurable partition $\alpha$ of $\mathbb{R}^d$,     its  \emph{entropy}   with respect to $(\varphi,\nu)$ is defined as 
\begin{align}
h_{\nu}(\varphi,\alpha):=\lim_{n\to+\infty}\frac{1}{n}\int_{\Omega}H_{\nu_{\omega}}\left(\bigvee_{i=0}^{n-1}\varphi(i,\omega)^{-1}\alpha\right)\,d\mathbb{P}(\omega),
\end{align}
where $H_{\nu_{\omega}}(\alpha)=\sum_{A\in\alpha}-\nu_{\omega}(A)\log\nu_{\omega}(A).$ 
\begin{definition}[Metric entropy of RDS]
\label{25-11-11-1122}
    The \emph{metric entropy} of  a measure-preserving RDS $(\varphi,\nu)$ on $\mathbb{R}^d$ is defined as
    $$h_{\nu}(\varphi):=\sup_{\alpha}h_{\nu}(\varphi,\alpha),$$
   where  $\alpha$ runs over finite Borel measurable partitions of $\mathbb{R}^d$. 
\end{definition} 

\begin{remark}
Throughout this paper, we mainly consider the time-1 map of the RDS generated by equation  \eqref{26-6-8-1611}.  In fact,  similar concepts also can be considered by  parallel argument. The reader can see the details in \cite[Chapter V]{MR1369243}.
\end{remark}


 We now recall the definition of an SRB measure in the random setting.

\begin{definition}
An  invariant measure $\nu\in\mathcal{I}(\varphi)$ is called  an SRB measure if the following conditions  hold:
\begin{enumerate}[(i)]
\item  the top Lyapunov exponents  $\chi^{1,+}_{\nu}:=\int_{\Omega\times\mathbb{R}^d}\chi^1_{\nu}(\omega,x)^+\,d\nu(\omega,x)>0$;
\item for every measurable partition $\eta$ subordinated to the unstable manifolds, its conditional measure $\nu^{\eta}$ is absolutely continuous with respect to  the volume measure induced by the unstable manifolds. 
\end{enumerate}
\end{definition}

  For the  RDS on $\mathbb{R}^d$, the entropy formula was also established    in \cite{MR3175262} under additional integrability estimates. By adjusting the arguments in \cite{MR1369243}, we can obtain the following characterization of the SRB measures in  non-compact phases. The needed  logarithmic integrability conditions are  given in \Cref{26-6-9-1713}.
\begin{thm}
\label{26-7-31-1029}
Let $\varphi$ be the RDS on $\mathbb{R}^d$ generated by \eqref{26-6-8-1611}.
\begin{itemize}
\item[(a)] An  invariant measure $\nu\in\mathcal{I}(\varphi)$  is an SRB measure  if and only if  $\chi^{1,+}_{\nu}>0$ and the entropy formula holds:
\begin{align}
\label{26-7-31-1028}
h_{\nu}(\varphi)=\int_{\Omega\times \mathbb{R}^d}\sum_{i=1}^d\chi^i_{\nu}(\omega,x)^+\, d\nu(\omega,x),
\end{align}
where $\chi^i_{\nu}$ is the Lyapunov exponents of $(\varphi,\nu)$ counted with multiplicities
\item[(b)] Let $\mu$ be the unique Markovian measure of \eqref{26-6-8-1611}.  If $\chi^{1,+}_{\mu}>0$, then $\mu$ is an SRB measure.
\end{itemize}
\end{thm}

\subsection{Main results}

We are now ready to state our main results.
\begin{thm}
\label{26-7-31-1031}
Assume that equation  \eqref{26-6-8-1611}  satisfies \hyperref[26-6-10-1501]{(H1)}-\hyperref[26-6-10-1503]{(H3)}.  Then the Markovian measure $\mu$ is the \textbf{unique invariant measure} which satisfies  the entropy formula \eqref{26-7-31-1028}.
\end{thm}

In the present setting,  \Cref{26-7-31-1031} relies on the following two properties.  

\noindent\textbf{Exponential ergodicity.}
Exponential ergodicity allows us to establish the multiplicative mean ergodic theorem in
\Cref{26-8-17-1609}. 

\noindent\textbf{Logarithmic integrability.}
The logarithmic integrability estimates in
\Cref{26-6-9-1713} allow us to apply Pesin theory and establish the volume  bounds of  Bowen ball needed to relate asymptotic log-moment
generating function, metric entropy, and Lyapunov exponents.
\begin{remark}
It is worth mentioning that the hypotheses of \Cref{26-7-31-1031} are only  imposed on the equation \eqref{26-6-8-1611}. No additional assumptions of hyperbolicity or positive entropy are required for the associated RDS.
\end{remark}
\begin{remark}
The assumptions
\hyperref[26-6-10-1501]{\textnormal{(H1)}}--\hyperref[26-6-10-1503]{\textnormal{(H3)}} are not intended to be optimal for establishing
the uniqueness of SRB measures.  For instance, the method can accommodate degenerate noise and coefficients of finite differentiability. Nevertheless, they provide a convenient framework that is sufficiently general for the applications considered in this paper. 
\end{remark}

Combining \Cref{26-7-31-1031}  with  \Cref{26-7-31-1029}, one immediately has following corollary.

\begin{cor}
\label{26-8-25-1504}
Under the setting of  \Cref{26-7-31-1031},    equation \eqref{26-6-8-1611} has an SRB measure if and only if $\chi^1_{\mu}>0$.  In particular,  equation  \eqref{26-6-8-1611} admits no SRB measures when $\chi^1_{\mu}\leq 0$.
\end{cor}
Now, we discuss  the following three cases of top sample Lyapunov exponents  and their connection with the existence results about SRB measures:

{\bf Case 1 $\chi^1_{\mu}>0$:}   $\mu$ is an SRB measure of  equation  \eqref{26-6-8-1611} following from entropy formula and its characterization (see \Cref{26-7-31-1029}).

{\bf Case 2 $\chi^1_{\mu}<0$:}  the compact random  attractor $\mathfrak{A}$ is a single random point under additional  conditions by synchronization \cite{MR696907,MR1144097,MR3630300,MR3663624,MR4791422}. It follows that   equation \eqref{26-6-8-1611} has no SRB measures. However, \Cref{26-8-25-1504} provides a more general condition to illustrate that  equation \eqref{26-6-8-1611} has no SRB measures.

{\bf Case 3 $\chi^1_{\mu}=0$:} this is a critical case. There are few literatures to describe the dynamical behaviour.  Our theory  implies that  equation  \eqref{26-6-8-1611}  also has no SRB measures under general conditions.

This paper is organized as follows.  In \Cref{26-8-29-2212}, we establish the integrability estimates required
for the RDS  generated by  equation \eqref{26-6-8-1611} and give  its Lyapunov spectrum.  In \Cref{26-8-17-1609}, we  develop the multiplicative mean ergodic theorem for  Feynman-Kac semigroup.  In \Cref{26-8-29-2220},  we establish a weak Katok's formula for non-ergodic invariant measure on a non-compact phase space  using the entropy theory and ergodic theory.  In \Cref{26-8-29-2229},  combining    non-uniform Pesin theory, Katok's formula and multiplicative mean ergodic theorem,  we complete the proof of   \Cref{26-7-31-1031}. Finally,  we apply
 \Cref{26-8-25-1504} to  obtain a dichotomy of SRB measures  for the 2D additive-noise stochastic Hopf normal form in \Cref{26-8-29-2233}.

\section{Multiplicative ergodic theorem and Lyapunov exponents}
\label{26-8-29-2212}
In this section, we verify the uniform log integrability estimate for the RDS $\varphi$ generated by equation \eqref{26-6-8-1611}. Consequently, we
can define the Lyapunov exponents for every invariant measure. 
\begin{proposition}
\label{26-6-9-1713}
The RDS $\varphi$ generated by \eqref{26-6-8-1611} admits a compact pullback random attractor $\mathfrak A$.
Moreover,  we have 
\begin{align}
\label{26-6-12-1605}
&\mathbb E\bigg[
 \sup_{0\leq t\leq1}
 \log^+\sup_{x\in\mathfrak A^1(\omega)}
 \|\mathrm d_x\varphi(t,\omega)\|
 +\sup_{0\leq t\leq1}
 \log^+\sup_{x\in\mathfrak A^1(\omega)}
 \|\nabla_x\mathrm d\varphi(t,\omega)\|
 \bigg]<\infty,\\
\label{26-9-4-2050}
&\mathbb E\bigg[
 \sup_{0\leq t\leq1}
 \log^+\sup_{x\in\mathfrak A^1(\omega)}
 \|\big(\mathrm d_x\varphi(t,\omega)\big)^{-1}\|
 \bigg]<\infty,
\end{align}
where  $\mathfrak A^1(\omega):=\{y\in\mathbb R^d:\dist(y,\mathfrak A(\omega))\leq1\}$.

\end{proposition}

\begin{proof}
For convenience, we write
\[
 X_t^x:=\varphi(t,\omega,x),\qquad
 J_t^x:=\mathrm d_x\varphi(t,\omega),\qquad
 Y_t^x:=(J_t^x)^{-1},\qquad
 K_t^x:=\nabla_x\mathrm d\varphi(t,\omega).
\]
Throughout the proof, $C,c>0$ denote deterministic constants which may change from line to line. Constants with subscripts may additionally depend on the indicated parameters. 

\begin{claim}
\label{claim:localized-C2-estimate}
For $r\geq1$, let
\begin{equation}
\label{eq:def-Zn-clean}
 Z_r:=1+\sup_{\substack{0\leq t\leq1\\|x|\leq r}}
 \big(\|J_t^x\|+\|Y_t^x\|+\|K_t^x\|\big).
\end{equation}
For every $q\geq1$, there exists $C_q>0$ such that  $\mathbb{E}[Z_r^q]\leq C_qe^{C_qr^2}$  whenever $r\geq1$.
\end{claim}

\begin{proof}
 Let $\mathcal L$ be the generator of equation \eqref{26-6-8-1611},  and  denote $V_\eta(x)=e^{\eta|x|^2}$ for  $\eta>0$. Since the diffusion
coefficients are bounded and $p>2$, we have
\begin{align*}
 \frac{\mathcal LV_\eta(x)}{V_\eta(x)} &=2\eta\langle x,b(x)\rangle+\eta\sum_{i=1}^m|\sigma_i(x)|^2  +2\eta^2\sum_{i=1}^m\langle x,\sigma_i(x)\rangle^2
\\ &\leq-\alpha\eta|x|^p+C_\eta.
\end{align*}
For any $\lambda>0$, choose $\eta$ such that $\alpha\eta\geq\lambda$.
Applying It\^o's formula to
\[
 \exp\left(\eta|X_t^x|^2+
 \lambda\int_0^t(1+|X_s^x|^p)\,ds-(C_\eta+\lambda)t\right)
\]
and using localization, we obtain
\begin{equation}
\label{eq:exponential-energy-section2}
 \mathbb{E}\left[\exp\left(\lambda\int_0^1(1+|X_s^x|^p)\,ds\right)\right]\leq C_\lambda e^{C_\lambda|x|^2}.
\end{equation}
The first variation and its inverse satisfy
\begin{align*}
 dJ_t^x
 &=\nabla b(X_t^x)J_t^x\,dt
   +\sum_{i=1}^m\nabla\sigma_i(X_t^x)J_t^x\,dW_t^i,\\
 dY_t^x&=Y_t^x\Big(-\nabla b(X_t^x)
       +\sum_{i=1}^m\big(\nabla\sigma_i(X_t^x)\big)^2\Big)dt
   -\sum_{i=1}^mY_t^x\nabla\sigma_i(X_t^x)\,dW_t^i,
\end{align*}
with $J_0^x=Y_0^x=\Id$. Applying It\^o's formula to $\log(1+\|U_t^x\|^2)$, with
$U=J,Y$, we have 
\[
 \log(1+\|U_t^x\|^2)
 \leq C+C\int_0^1(1+|X_s^x|^p)\,ds+M_t^{U,x},
 \qquad 0\leq t\leq1,
\]
where $M^{U,x}$ is a continuous local martingale starting at zero
with $\langle M^{U,x}\rangle_1\leq C$.   For any $q\geq1$,  the exponential martingale maximal inequality, H\"older's inequality and \eqref{eq:exponential-energy-section2} therefore yield
\begin{equation}
\label{eq:JY-moment-clean}
 \mathbb E\big[\sup_{0\leq t\leq1}\big(\|J_t^x\|^q+\|Y_t^x\|^q\big) \big]\leq C_qe^{C_q|x|^2}.
\end{equation}

For $u,v\in\mathbb R^d$, the variation-of-constants formula for
the second variation gives
\begin{align}
 \notag K_t^x(u,v)=J_t^x\bigg(&\int_0^tY_s^x\big(
 \nabla^2b(X_s^x)(J_s^xu,J_s^xv)
 -\sum_{i=1}^m\nabla\sigma_i(X_s^x)
       \nabla^2\sigma_i(X_s^x)(J_s^xu,J_s^xv)
 \big)\,ds\\
 &\label{26-9-19-2159}+\sum_{i=1}^m\int_0^tY_s^x
       \nabla^2\sigma_i(X_s^x)(J_s^xu,J_s^xv)\,dW_s^i
 \bigg).
\end{align}
By \hyperref[26-6-10-1502]{(H2)}, the time integral is bounded by
\[
 C|u||v|\left(\int_0^1(1+|X_s^x|^p)\,ds\right)
 \left(\sup_{0\leq s\leq1}\|J_s^x\|\right)^3
 \sup_{0\leq s\leq1}\|Y_s^x\|.
\]
The Burkholder--Davis--Gundy inequality controls the stochastic
integral. Together with H\"older's inequality,
\eqref{eq:exponential-energy-section2} and
\eqref{eq:JY-moment-clean}, we can obtain the  same moment estimate for $K_t^x$. Moreover, It\^o's formula for
$(1+|X_t^x|^2)^{q/2}$ and the Burkholder--Davis--Gundy inequality imply that 
$$\mathbb E\big[\sup_{0\leq t\leq1}|X_t^x|^q\big]\leq C_q(1+|x|^q),\quad  q\geq 1.$$
For all, we have 
\begin{equation}
\label{eq:pointwise-variations-section2}
 \mathbb E\Big[\sup_{0\leq t\leq1}
 \big(|X_t^x|^q+\|J_t^x\|^q+\|Y_t^x\|^q+\|K_t^x\|^q\big)\Big]\leq C_qe^{C_q|x|^2},\qquad q\geq1.
\end{equation}

We next obtain estimates uniform over deterministic balls. Fix $\ell>\max\{d,q,2\}$. For each $U\in\{X,J,Y\}$, the Sobolev inequality on unit balls
and Jensen's inequality yield
\begin{align}
\label{eq:local-sobolev-section2}
\mathbb E\Big[\sup_{\substack{0\leq t\leq1\\|x|\leq r}}|U_t(x)|^q \Big]\leq C_{q,\ell}\left(\int_{B(0,r+1)}\mathbb{E}\Big[\sup_{0\leq t\leq1}\bigl(|U_t(x)|^\ell+\|\nabla_xU_t(x)\|^\ell\bigr)\Big]\,dx\right)^{q/\ell},
\end{align}
where $r\geq1$. The pointwise moment estimates in
\eqref{eq:pointwise-variations-section2}, together with H\"older's inequality, also control the spatial derivatives. Thus,
\[
\mathbb {E}\Big[\sup_{0\leq t\leq1}\bigl(|U_t(x)|^\ell+\|\nabla_xU_t(x)\|^\ell\bigr)\Big]\leq C_\ell e^{C_\ell|x|^2},
\qquad U\in\{X,J,Y\}.
\]
It follows that 
\begin{equation}
\label{eq:uniform-first-variations-section2}
\mathbb E\Big[\sup_{\substack{0\leq t\leq1\\|x|\leq r}}\bigl(|X_t^x|^q+\|J_t^x\|^q+\|Y_t^x\|^q\bigr)\Big]
\leq C_qe^{C_qr^2},
\qquad r\geq1.
\end{equation}

To estimate $K$ uniformly in $x$, we use
\eqref{26-9-19-2159} again. The time integral is controlled by
\eqref{eq:uniform-first-variations-section2} and H\"older's
inequality. For the stochastic integral, fix $u,v$ among the
coordinate vectors and write
\[
 H_{i,s}^x:=Y_s^x\nabla^2\sigma_i(X_s^x)(J_s^xu,J_s^xv),
 \qquad
 N_t^x:=\sum_{i=1}^m\int_0^tH_{i,s}^x\,dW_s^i.
\]
The spatial derivative of $H_{i,s}^x$ is a sum of products involving
only $J_s^x,Y_s^x,K_s^x$ and the bounded derivatives
$\nabla^2\sigma_i,\nabla^3\sigma_i$. Stochastic Fubini and integration by parts against compactly supported smooth test functions imply that
\[
 \nabla_xN_t^x=\sum_{i=1}^m
 \int_0^t\nabla_xH_{i,s}^x\,dW_s^i
\]
in the weak sense. These operations are justified on every deterministic ball by \eqref{eq:pointwise-variations-section2}.
Thus, the Burkholder--Davis--Gundy inequality yields
\[
 \mathbb E\Big[\sup_{0\leq t\leq1}\big(|N_t^x|^\ell+\|\nabla_xN_t^x\|^\ell\big)\Big] \leq C_\ell e^{C_\ell|x|^2}.
\]
Applying \eqref{eq:local-sobolev-section2} to $N$, and  using \eqref{eq:uniform-first-variations-section2} for the factor $J_t^x$ in \eqref{26-9-19-2159},  we obtain \Cref{claim:localized-C2-estimate}.
\end{proof}

\begin{claim}
\label{25-6-12-1553}
The RDS $\varphi$ admits a compact random attractor $\mathfrak A$. Moreover, there exists a random variable $\mathsf R:\Omega\to[1,\infty)$ such that $\mathfrak A^1(\omega)\subset B(0,\mathsf R(\omega))$
for $\mathbb P$-a.s. $\omega\in\Omega$, and, for every $q\geq1$,
\begin{equation}
\label{eq:absorbing-radius-moment}
 \mathbb E[\mathsf R^q]<\infty.
\end{equation}
\end{claim}
\begin{proof}
Let $U_t^x=(1+|X_t^x|^2)^{1/2}$. According to  It\^o's formula and
\hyperref[26-6-10-1502]{(H2)}, one has 
\begin{equation}
\label{eq:Y-dissipative-comparison}
 dU_t^x\leq\big(C-c(U_t^x)^{p-1}\big)\,dt+dM_t^x,
 \qquad d\langle M^x\rangle_t\leq C\,dt,
\end{equation}
where $M^x$ is a continuous local martingale with $M_0^x=0$. For sufficiently large $r$, put $\delta_r=\kappa r^{-(p-2)}\leq1$, with $\kappa>0$ to be chosen. The process $z_t:=U_t^x-M_t^x$ is absolutely continuous.
On the event
\[
\left\{\sup_{0\leq t\leq\delta_r}|M_t^x|\leq r/16\right\},
\]
whenever $z_t\geq r/4$, we have
\[
U_t^x=z_t+M_t^x
\geq z_t-r/16\geq3z_t/4.
\]
Consequently, for sufficiently large $r$ we have
\[
z_t-z_s\leq-c'\int_s^t z_{\tau}^{p-1}\,d\tau
\]
whenever $0\leq s\leq t\leq\delta_r$ and $z_{\tau}\geq r/4$ throughout $[s,t]$. The scalar comparison estimate shows that the time needed
to reach $r/4$ from any larger initial value is at most
\[
\frac{(r/4)^{-(p-2)}}{c'(p-2)},
\]
and $z$ cannot subsequently cross above $r/4$. By choosing
$\kappa$ sufficiently large, we have 
\[
|X^x_{\delta_r}|\leq U_{\delta_r}^x=z_{\delta_r}+M_{\delta_r}^x\leq r/4+r/16=5r/16<r/2
\]
on this event. Since $d\langle M^x\rangle_t\leq C\,dt$, the exponential martingale maximal inequality gives
\begin{align}
\label{eq:pointwise-entrance-short}
\sup_{|x|\leq2r}\mathbb P\{|X_{\delta_r}^x|\geq r/2\}&\leq\sup_{|x|\leq2r}\mathbb P\left\{\sup_{0\leq t\leq\delta_r}|M_t^x|>r/16\right\} 
\\\notag&\leq Ce^{-cr^2/\delta_r}
\\\notag&\leq Ce^{-cr^p}.
\end{align}
By \Cref{claim:localized-C2-estimate}, we have $\mathbb E[Z_{2r}]\leq Ce^{C_1r^2}$. Choose $C_2>C_1+1$, and take an $e^{-C_2r^2}$-net of $B(0,2r)$ containing at most
$Cr^de^{dC_2r^2}$ points.
If $Z_{2r}\leq re^{C_2r^2}/2$ and the image of every net point at time $\delta_r$ belongs to $B(0,r/2)$, the mean-value theorem shows that
the image of the whole ball belongs to $B(0,r)$.  Markov's inequality and \eqref{eq:pointwise-entrance-short} therefore give
\begin{align}
\label{eq:uniform-entrance-short}
 \mathbb P\{\varphi(\delta_r,\omega,B(0,2r))\not\subset B(0,r)\}\leq\frac{2e^{-C_2r^2}}{r}\mathbb E[Z_{2r}]+Cr^de^{dC_2r^2-cr^p}\leq Ce^{-c_0r^2},
\end{align}
where we used $p>2$.  

We now construct a compact random absorbing set.  For each sufficiently large integer $r$, set
\[
 r_j:=2^jr,\qquad \delta_{r_j}
       =\kappa r_j^{-(p-2)},\qquad
 S_j:=\sum_{i=0}^j\delta_{r_i},\qquad S_{-1}:=0,
\]
and
\[
 T_r:=\sum_{j=0}^{\infty}\delta_{r_j}=\frac{\kappa r^{-(p-2)}}{1-2^{-(p-2)}}<\infty.
\]
Define
\[
 G_r:=
 \bigcap_{j=0}^{\infty}
 \left\{
 \varphi\bigl(\delta_{r_j},\theta_{-S_j}\omega,
                  \overline B(0,2r_j)\bigr)
 \subseteq\overline B(0,r_j)
 \right\}.
\]
By \eqref{eq:uniform-entrance-short} and invariance of $\mathbb{P}$ with respect to $\theta$, we have
\[
 \mathbb P(G_r^c)
 \leq C\sum_{j=0}^{\infty}e^{-c_0r_j^2}
 \leq C'e^{-c'r^2}.
\]
Fix $\omega\in G_r$ and $x\in\mathbb R^d$.  It is clear that 
\[
 x_N:=\varphi(T_r-S_N,\theta_{-T_r}\omega,x)\to x\quad\text{as } N\to\infty.
\]
Thus $x_N\in\overline B(0,2r_N)$ for all sufficiently large $N$.  Applying the cocycle property successively for $j=N,N-1,\ldots,0$, we obtain 
\[
 \varphi(T_r,\theta_{-T_r}\omega,x)\in\overline B(0,r).
\]
By continuity of $\varphi(t,\cdot)$,  for any $\omega\in G_r$ we have  $\varphi(T_r,\theta_{-T_r}\omega,\mathbb R^d)
 \subseteq\overline B(0,r).$ It follows that for every $t\geq T_r$
 \begin{align*}
 \varphi(t,\theta_{-t}\omega,\mathbb R^d)&=\varphi\bigl(T_r,\theta_{-T_r}\omega,\varphi(t-T_r,\theta_{-t}\omega,\mathbb R^d)\bigr)\subseteq\overline B(0,r).
 \end{align*}
Choose an integer $n_0$ sufficiently large that all the preceding estimates hold for $r\geq n_0$ and $T_{n_0}\leq1$.
Define
\[
 \mathsf R_0(\omega)
 :=\min\{r\in\mathbb N:r\geq n_0,\ \omega\in G_r\}.
\]
Since  for every integer $r\geq n_0$,
\[
 \mathbb P\{\mathsf R_0>r\}
 \leq\mathbb P(G_r^c)\leq C'e^{-c'r^2},
\]
we have $\mathbb E[\mathsf R_0^q]<\infty$ for every $q\geq1$. For $\mathbb{P}$-a.s. $\omega\in\Omega$, we have
\[
 \varphi(t,\theta_{-t}\omega,\mathbb R^d)\subseteq\overline B(0,\mathsf R_0(\omega))
,\qquad t\geq1.
\]
It therefore pullback absorbs every bounded non-random set. By \cite[Theorem~3.11]{MR1305587},  the RDS $\varphi$ admits a compact random attractor $\mathfrak{A}$. By strict invariance and the preceding inclusion at time $1$, one has that for $\mathbb{P}$-a.s. $\omega\in\Omega$
\[
 \mathfrak A(\omega)=\varphi\bigl(1,\theta_{-1}\omega, \mathfrak A(\theta_{-1}\omega)\bigr) \subseteq\overline{B}(0,\mathsf R_0(\omega)).
\]
Hence we obtain \Cref{25-6-12-1553} by letting $\mathsf R:=\mathsf R_0+2$.
\end{proof}

We now complete the proof of \Cref{26-6-9-1713}. By \Cref{claim:localized-C2-estimate}, there is $C_3>0$
such that $\mathbb E[Z_n]\leq Ce^{C_3n^2}$ for $n\geq1$.
Choose $C_4>C_3$ and define
\[
 S:=1+\sum_{n=1}^\infty e^{-C_4n^2}Z_n.
\]
Then $\mathbb E[S]<\infty$, and $Z_r\leq Z_{\lceil r\rceil}\leq e^{C_4\lceil r\rceil^2}S$ for any $r\geq1$. In particular, this pathwise bound applies to the random radius
$\mathsf R$, without any independence assumption. Hence
\[
 \mathbb{E}[\log Z_{\mathsf R}]\leq C_4\mathbb{E}\big[\lceil\mathsf R\rceil^2\big]+\mathbb{E}[\log S]
 \leq C_4\mathbb{E}\big[\lceil\mathsf R\rceil^2\big]+\log\mathbb E[S]
 <\infty.
\]
Since $\mathfrak A^1(\omega)\subset B(0,\mathsf R(\omega))$,
each of the three logarithmic terms in
\eqref{26-6-12-1605} and \eqref{26-9-4-2050} is bounded by $\log Z_{\mathsf R}$. This proves  \Cref{26-6-9-1713}.
\end{proof}

The uniform integrability established in \Cref{26-6-9-1713}
allows us to apply the multiplicative ergodic theorem to the
derivative cocycle with respect to every invariant measure.
Indeed, every $\nu\in\mathcal I(\varphi)$ is supported on the
random attractor, that is, $\nu(\mathfrak A)=1$.
To see this, for every $R,\epsilon>0$ and
$t\geq0$, one has 
\begin{align*}
 &\nu\{(\omega,x):\dist(x,\mathfrak A(\omega))>\epsilon\}\\
 &\qquad\leq\nu\{(\omega,x):|x|>R\}
 +\mathbb P\{\dist(\varphi(t,\theta_{-t}\omega,B(0,R)),
                         \mathfrak A(\omega))>\epsilon\}.
\end{align*}
Letting first $t\to\infty$ and then $R\to\infty$ proves the assertion. Since the marginal of $\nu$ is $\mathbb P$, it follows that
\[
 \mathbb E_\nu\left[\sup_{0\leq t\leq1}
 \left(\log^+\|\mathrm d_x\varphi(t,\omega)\|
 +\log^+\|(\mathrm d_x\varphi(t,\omega))^{-1}\|\right)\right]
 <\infty.
\]
Therefore,  we can  apply the multiplicative ergodic theorem (see, for example, \cite{MR2674952,MR1369243,A}) in the following form.

\begin{thm}
\label{26-6-13-1113}
For any $\varphi$-invariant measure $\nu$, there  is  a $\Phi$-invariant Borel subset  $\Gamma$ of $\Omega\times \mathbb{R}^d$  with $\nu(\Gamma)=1$,  a measurable function  $\mathbf{r}:\Gamma\to\mathbb{N}$,   a collection of measurable functions $\lambda_{\nu}^1,\dots,\lambda_{\nu}^{\mathbf{r}}: \Gamma\to\mathbb{R}$ with  
$$\lambda_{\nu}^1(\omega,x)>\lambda_{\nu}^2(\omega,x)>\cdots>\lambda_{\nu}^{\mathbf{r}}(\omega,x),$$
and  a measurable decomposition $\mathbb{R}^d=E^{1}(\omega,x)\oplus\cdots\oplus E^{\mathbf{r}}(\omega,x)$  with the following properties 
  \begin{enumerate}[\rm (I)]
  \item  Invariance:  for each $j=1,2,\dots, \mathbf{r}(\omega, x)$, one has that   
  $$\mathrm{d}_x\varphi(t,\omega)E^{j}(\omega,x )=E^{j}\big(\theta_t\omega, \varphi(t,\omega)x\big);$$
  \item Lyapunov exponent: 
  $$\lim_{t\to\infty}\frac{\log|\mathrm{d}_x\varphi(t,\omega)v|}{t}=\lambda^j_{\nu}(\omega,x)$$ 
  for any  $ v\in E^j(\omega,x)\setminus\{0\}$;
     \item   Temperedness: the projection operators associated with the decompositions 
  \begin{align*}
  \mathbb{R}^d=\big(\bigoplus_{i=1}^jE^i(\omega,x)\big)\oplus\big(\bigoplus_{i=j+1}^{\mathbf{r}}E^i(\omega,x)\big)
  \end{align*} 
  are  strongly measurable and tempered.
  \end{enumerate}
  \end{thm}
\begin{remark}
Throughout this paper,  we always  denote  the Lyapunov exponents of measure-preserving RDS $(\varphi,\nu)$ counted with multiplicity as $\chi^1_{\nu}\geq \cdots\geq\chi_{\nu}^{d}$. If we additionally assume that $\nu$ is ergodic, then the Lyapunov exponents and dimensions of the Oseledets spaces are constants up to a $\nu$-null measure set. 
\end{remark}

\section{Multiplicative mean ergodic theorem}
\label{26-8-17-1609}
Kontoyiannis and Meyn \cite{MR1952001,KontoyiannisMeyn2005} established multiplicative mean ergodic theorems of Markov processes  
through the spectral analysis of Feynman--Kac operators. In this section, we adapt this approach to prove a  multiplicative mean ergodic theorem for   potentials depending on both the state and a finite future interval of the driving noise, using the weighted Harris theorem \cite{MeynTweedie2009,MR2857021} and analytic perturbation theory \cite[Chapter~VII]{Kato1995}. Throughout this paper,  a  function $f:\Omega\times\mathbb R^d\to\mathbb R$ is called a \emph{potential function} if it satisfies:
\begin{itemize}
\item $f$ is finite-future measurable,  i.e.  there is some $\tau\geq0$ such that $f$ is $\mathscr{F}_0^{\tau}\otimes\mathscr B(\mathbb R^d)$-measurable;
\item for $\mathbb P$-a.s. $\omega\in\Omega$, the function
$f_\omega:=f(\omega,\cdot)$ is bounded and Lipschitz, and
\[
\|f\|:=\operatorname*{ess\,sup}_{\omega\in\Omega}
\bigl(\|f_\omega\|_\infty+\operatorname{Lip}(f_\omega)\bigr)<\infty.
\]
\end{itemize}
Our multiplicative mean ergodic theorem is  stated as following: 
\begin{thm}
\label{prop:ff-exponential-bound}
Let $f:\Omega\times\mathbb R^d\to\mathbb R$ be a potential function. Then there exists $s_0>0$ such that for any sufficiently large $\tau>0$  there is  a real-analytic function $\Lambda_f: (-s_0,s_0)\to\mathbb R$ such that:
\begin{enumerate}[(a)]
\item $\Lambda_f'(0)=\mathbb{E}_\mu[f]$, where $\mu$ is the Markovian measure of the equation \eqref{26-6-8-1611};
\item for every $|s|<s_0$ and $x\in\mathbb R^d$, one has that
\begin{equation}\label{26-9-30-0951}
\lim_{n\to\infty}\frac1n\log\mathbb E\left[
\exp\left(s\sum_{i=0}^{n-1}f(\Phi_{i\tau}(\omega,x))\right)\right]=\Lambda_f(s),
\end{equation}
and the convergence is locally uniform in $x\in\mathbb R^d$.
\end{enumerate}
\end{thm}

\begin{proof}
Recall that $\gamma$ is the unique stationary measure of
\eqref{26-6-8-1611}, and let $(\mathcal Q_t)_{t\geq0}$ denote its
Markov semigroup. 

\begin{claim}\label{clm:weighted-mixing}
Let $V(x)=1+|x|^2$. Then $\gamma(V)<\infty$, and there exist constants $\lambda, C>0$  such that
\[
\|\mathcal Q_tu-\gamma(u)\mathbf1\|_V
\leq C e^{-t\lambda}\|u\|_V,
\]
for any $t\geq0$ and  every complex-valued Borel function $u$ with
\[
\|u\|_V:=\sup_{x\in\mathbb R^d}\frac{|u(x)|}{V(x)}<\infty.
\]
\end{claim}
\begin{proof}
Dissipativity and boundedness of the diffusion coefficients give
\[
\mathcal LV(x)
=2\langle x,b(x)\rangle+\sum_{i=1}^m|\sigma_i(x)|^2
\leq-2\alpha|x|^p+C_0\leq-c_1V(x)+C_1
\]
for some $c_1,C_1>0$. By  Dynkin's formula, there exist constants $C>0$ and $0<c<1$ such that 
\[
\mathcal Q_1V\leq cV+C.
\]
Smoothness and uniform ellipticity imply that the time-one transition density  is continuous and strictly positive.  It follows the existence of  a small set for $\mathcal Q_1$. Positivity also gives irreducibility and aperiodicity. We obtain \Cref{clm:weighted-mixing} immediately by the weighted Harris theorem \cite[Chapters~15--16]{MeynTweedie2009}. 
\end{proof}

Define the weighted Banach space
\[
 \mathcal{B}_V:=\left\{g:\mathbb{R}^d\to\mathbb C\ \text{is measurable}: \|g\|_V<\infty\right\},
\]
endowed with the  weighted norm
\[
 \|g\|_V :=\sup_{x\in\mathbb{R}^d}\frac{|g(x)|}{V(x)}.
\]
Without loss of generality,  we assume that $\tau$ is  a positive constant such that $f$ is $\mathscr{F}_0^{\tau}\otimes\mathscr{B}(\mathbb{R}^d)$-measurable.  For convenience,  for any $t\geq 0, x\in\mathbb{R}^d$ denote $X_t^x=\varphi(t, \cdot, x)$, and 
\[
S_nf:=\sum_{i=0}^{n-1}f\circ\Phi_{i\tau}.
\]
For $s\in\mathbb{R}$, consider the following  Feynman--Kac operator  $\mathcal Q^{(s)}u(x):=\mathbb{E}[e^{sf(\omega,x)}u(X_{\tau}^x)]$, where $u\in \mathcal{B}_V$. By independence of Brownian increments and the cocycle property of the RDS $\varphi$,  we have
\begin{equation}\label{eq:ff-one-step-fk}
\mathcal Q_n^{(s)}u(x)=\mathbb E\left[e^{sS_nf(\omega,x)}u(X_{n\tau}^x)\right],
\end{equation}
where  $\mathcal Q_n^{(s)}$ is  the $n$-th composition of $\mathcal Q^{(s)}$.  For any $u\in\mathcal{B}_V$ and $s\in\mathbb{R}$, one has that 
\begin{align*}
\Big|\sum_{k=0}^{\infty}\frac{s^k}{k!}\mathbb{E}\left[f(\omega,x)^k u(X_{\tau}^x)\right]\Big|\leq &\sum_{k=0}^{\infty}\frac{|s|^k}{k!}\mathbb{E}\left[|f(\omega,x)|^k\cdot |u(X_{\tau}^x)|\right]
\\\leq&\sum_{k=0}^{\infty}\frac{|s|^k}{k!}\mathcal{Q}_{\tau}V(x)\cdot \|f\|^k\cdot \|u\|_{V}
\\<&\infty.
\end{align*}
Thus we obtain
\[
\mathcal Q^{(s)}u(x)=\sum_{k=0}^{\infty}\frac{s^k}{k!}\mathbb E\left[f(\omega,x)^k u(X_{\tau}^x)\right].
\]
It  follows that  $s\mapsto \mathcal{Q}^{(s)}_{}$ is analytic in the operator norm with $\mathcal{Q}_1^{(0)}=\mathcal{Q}^{(0)}=\mathcal{Q}_{\tau}$. By \Cref{clm:weighted-mixing},  we know that $1$ is a simple isolated eigenvalue of $\mathcal Q_{\tau}$, whose spectral projection is
$$\mathcal{R} u=\gamma(u)\mathbf{1}.$$  Applying  the analytic perturbation theory (\cite[Chapter~VII]{Kato1995}), there exists  an analytic simple eigenvalue $\kappa(s)$ and an analytic
rank-one spectral projection $\mathcal R_s$ satisfying that 
\[
\kappa(0)=1,\quad\text{and}\quad \mathcal R_0=\mathcal{R}.
\]
Then there exists  $s_0>0$ such that  for any  $|s|<s_0$
\begin{itemize}
 \item the eigenvalue $\kappa(s)$ of $\mathcal{Q}^{(s)}$ is positive;
 \item  there exists constants  $\lambda_s, C_s>0$ such that for any $n\in\mathbb{N}$ one has 
\begin{equation}\label{eq:ff-one-step-spectral}
\left\|\kappa(s)^{-n}\mathcal Q_{n}^{(s)}-\mathcal R_s\right\|_V\leq C_se^{-n\lambda_s}.
\end{equation}

\end{itemize} It is clear that $h_s:=\mathcal{R}_s\mathbf{1}$  is an eigenfunction of $\mathcal{Q}^{(s)}$ with respect to $\kappa(s)$.
Differentiating $\mathcal Q^{(s)}h_s=\kappa(s)h_s$, one has that 
\[
\mathbb E[f(\omega,\cdot)]+\mathcal Q_{\tau}(h'_0)=\kappa'(0)+h'_0.
\]
Therefore,  we obtain $\kappa'(0)=\mathbb E_\mu[f]$,  since $\mu_{\omega}$ is independent of $f$ and $\mathbb{E}[\mu_{\cdot}]=\gamma$. This proves  (a) of \Cref{prop:ff-exponential-bound} by taking $\Lambda_f(s)=\log \kappa(s)$.

According to the positivity of $\mathcal Q^{(s)}$ and \eqref{eq:ff-one-step-spectral}, one has that for any $x\in\mathbb{R}^d$
\[
h_s(x)=\mathcal{R}_s\mathbf{1}=\lim_{n\to\infty}\kappa(s)^{-n}\mathcal Q_{n}^{(s)}\mathbf1(x)\geq0.
\]
Since  $\lim_{s\to 0}\gamma(h_s)=1$,  we can  decrease  $s_0$ such that  $\gamma(h_s)>0$ for any  $|s|<s_0$. The eigenfunction equation  yields for any $n\in\mathbb{N}$
\begin{align*}
\kappa(s)^nh_s&=\mathcal Q_{n}^{(s)}h_s\\
&\geq e^{-n|s|\cdot\|f\|_{\infty}}\mathcal{Q}_{n\tau}h_s\\
&\tag{By  \Cref{clm:weighted-mixing}}\geq e^{-n|s|\cdot\|f\|_{\infty}}\bigl(\gamma(h_s)-Ce^{-n\tau\lambda}\|h_s\|_V V\bigr).
\end{align*}
For every compact set $K\subset\mathbb R^d$,  we choose $N$ large enough such that for any $n\geq N$
\[
\inf_{x\in K}h_s(x)
\geq\kappa(s)^{-n}e^{-n|s|\cdot\|f\|_{\infty}}\frac{\gamma(h_s)}{2}>0.
\]
Applying \eqref{eq:ff-one-step-fk} and \eqref{eq:ff-one-step-spectral} to $u=\mathbf1$, we have
\[
\lim_{n\to\infty}\sup_{x\in K}
\left|\frac1n\log\mathbb E[e^{sS_nf(\omega,x)}]-\Lambda_f(s)\right|=0.
\]
This proves \eqref{26-9-30-0951} and its local uniformity in $x$.
\end{proof}

\section{Revisit of Katok's formula for RDS}
\label{26-8-29-2220}
 Katok's  formula  says that we can use the maximal cardinality of separated sets to approximate the metric entropy.  It was proved by Katok in \cite{MR573822}.   Subsequently, there is a lot of literature to investigate the Katok's  formula  (for example, see \cite{MR3834712, MR2511539, MR3485864, MR3062903}).  In this section, we   build a weak version of Katok's formula for  non-compact phase space $\mathbb{R}^d$ and non-ergodic invariant measure. The main methods rely on the ergodic decomposition   and construction of null-boundary partitions. Throughout this paper,  we always assume that $\varphi$ is the   RDS over    $(\Omega,\mathscr{F},\mathbb{P}, \theta)$ generated by \eqref{26-6-8-1611} and $\nu\in\mathcal{I}(\varphi)$.

To deal with the non-ergodic issue, we recall the ergodic decomposition of an invariant measure. The reader can refer to  the details in \cite[Theorem 6.2]{EW}, \cite[Theorem 2.33]{MET} and \cite[Theorem 5.4 in Chapter 0]{MR1369243}.  For any $\nu\in\mathcal{I}(\varphi)$, let 
$$\nu=\int_{\Omega\times\mathbb{R}^d}\nu^{\mathbf{x}}\,d\nu(\mathbf{x})$$
be the ergodic decomposition of $\nu$.  It means that 
\begin{itemize}
 \item  $\nu^{\mathbf{x}}$ is an ergodic measure of $\varphi$ and $\nu^{\Phi_1(\mathbf{x})}=(\Phi_1)_*\nu^{\mathbf{x}}=\nu^{\mathbf{x}}$ for $\nu$-a.s. $\mathbf{x}\in\Omega\times\mathbb{R}^d$;
  \item for any  $\Phi_1$-invariant measurable subset $A$, $\nu^{\mathbf{x}}(A)=1$ for $\nu$-a.s. $\mathbf{x}\in A$;
  \item  $h_{\nu}(\varphi)=\int_{\Omega\times\mathbb{R}^d}h_{\nu^{\mathbf{x}}}(\varphi)\,d\nu(\mathbf{x}).$
\end{itemize}  

Now, we recall  the definition of  Bowen  balls, which  is a typical object to capture the local dynamical behavior. Given $(\omega,x)\in\Omega\times\mathbb{R}^d$, $n\geq 0$ and $\epsilon>0$, the Bowen ball is defined as 
\begin{align}
B(\omega,x,n,\epsilon):=\{y\in\mathbb{R}^d:  |\varphi(i,\omega,y)-\varphi(i,\omega,x)|\leq\epsilon\,\text{ for }\, 0\leq i\leq n-1\}.
\end{align}
 For any subset $\Gamma$ of $\Omega\times\mathbb{R}^d$, let 
 $$\Gamma_{\omega}:=\{x\in\mathbb{R}^d: (\omega,x)\in\Gamma\}.$$
 A subset  $E$  of $\mathbb{R}^d$ is called an  $(\omega,n,\epsilon)$-separated set if  $y\notin B(\omega,x,n, \epsilon)$ for any $x\neq y\in E$.  Denote
 $$s(\omega,n,\epsilon;\Gamma_\omega):=\sup\{\#E:   E \text{ is a $(\omega,n,\epsilon)$-separated set and } E\subset \Gamma_{\omega}\}.$$ 
The main result of this section is stated as follows. 
\begin{thm}
\label{26-5-5-1704}
Let $\Gamma$ be a  positive $\nu$-measure subset.  Then there is a positive $\mathbb{P}$-measure  subset $\Omega_0$ of $\pi_{\Omega}(\Gamma)$ such that for any $\omega\in\Omega_0$
$$\lim_{\epsilon\to 0}\liminf_{n\to\infty}\frac{\log s(\omega,n,\epsilon;\Gamma_\omega)}{n}\geq \nu(\Gamma)^{-1}\int_\Gamma h_{\nu^{\mathbf x}}(\varphi)\,d\nu(\mathbf x).$$
\end{thm}
\begin{remark}
\label{26-9-5-1802}
If $\Gamma$ is a compact subset of $\Omega\times\mathbb{R}^d$,  then  $\omega\mapsto\log s(\omega,n,\epsilon;\Gamma_\omega)$ is $\mathscr{F}$-measurable  by \cite[Lemma 1.2.2]{KiL}.
\end{remark}

Before proving \Cref{26-5-5-1704}, we give some auxiliary lemmas.  For a finite measurable partition $\alpha$ of $\mathbb{R}^d$, denote its boundary by 
$$\partial\alpha=\bigcup_{A\in\alpha}\partial A.$$ 
This partition  is said to have a  \emph{$\gamma$-null boundary} if $\gamma(\partial \alpha)=0$ for some probability measure $\gamma$ on $\mathbb{R}^d$.  Combining Ruelle's inequality  \cite{MR3233706} or \cite[Proposition 4.4]{Lianzhang} with the  integrability conditions in \Cref{26-6-9-1713},   we know that $h_{\nu}(\varphi)<\infty$. The following lemma then follows  from  \cite[Lemma 4.1]{Lianzhang}.
\begin{lemma}
\label{26-2-10-2156}
  For any small $\delta>0$,  there exists a finite measurable partition $\alpha$ of $\mathbb{R}^d$ with $(\pi_{\mathbb{R}^d})_*\nu$-null boundary  such that $h_{\nu}(\varphi)\leq h_{\nu}\big(\varphi,\alpha\big)+\delta.$
 \end{lemma}

Fix a finite measurable partition $\alpha$ of $\mathbb{R}^d$ and $\omega\in\Omega$.  Denote
\[
 \alpha_\omega^n
 :=\bigvee_{i=0}^{n-1}
       \bigl(\varphi(i,\omega)\bigr)^{-1}\alpha, 
\]
and  let $\alpha_\omega^n(x)$ be the atom of $\alpha^n_{\omega}$ that  contains   $x$.  This measurable partition  induces a  family of pseudo-metrics on $\mathbb{R}^d$ as follows 
\begin{align}
\label{26-9-26-1004}
\rho_{\alpha^n_{\omega}}(x,y)=\frac{1}{n}\sum_{i=0}^{n-1}\rho_{\alpha}(\varphi(i,\omega)x,\varphi(i,\omega)y),
\end{align}
where $\rho_\alpha$ is the discrete pseudo-metric induced by $\alpha$, recording whether two points belong to the same atom of $\alpha$.  The following lemma  bounds the number of atoms of $\alpha_\omega^n$ that intersect a closed ball in the pseudo-metric \eqref{26-9-26-1004}.
\begin{lemma}[\cite{MR573822} or {\cite[Lemma A.1]{MR3233706}}]
\label{26-2-10-2333}
 Given a  small positive constant $r$,  for any $(\omega,x)\in\Omega\times\mathbb{R}^d$ and $n\in\mathbb{N}$ one has that 
$$\#\{ A\in\alpha^n_{\omega} :  A\cap\mathbf{B}_{\rho_{\alpha^n_{\omega}}}(x,r)\neq\varnothing\}\leq \mathrm{B}(r,|\alpha|, n):=\sum_{i=0}^{\lfloor nr\rfloor}\big(|\alpha|-1\big)^i\binom{n}{i},$$
where $\mathbf{B}_{\rho_{\alpha^n_{\omega}}}(x, r)$ is  the closed ball in the pseudo-metric $\rho_{\alpha^n_{\omega}}$ centered at $x$ with radius $r$.
\end{lemma}


After completing all preparations, we give the final proof of \Cref{26-5-5-1704}.
\begin{proof}[Proof of \Cref{26-5-5-1704}]

Firstly, we assume that there is a positive constant $h_*$ such that $h_{\nu^{\mathbf{x}}}(\varphi)\geq h_*$ for $\nu$-a.s. $\mathbf{x}\in\Gamma$.  Let
\[
  \mathcal C:=\left\{\mathbf x\in\Omega\times\mathbb R^d:
       h_{\nu^{\mathbf x}}(\varphi)<h_*\right\}.
\]
It follows that  $\nu(\Gamma\cap\mathcal C)=0.$   Since $\nu^{\Phi(\mathbf x)}=\nu^{\mathbf x}$ for $\nu$-a.s. $\mathbf{x}\in\Omega\times\mathbb{R}^d$, the set $\mathcal C$ is $\Phi_1$-invariant modulo a $\nu$-null set.  Then we obtain
\begin{align*}
  0<\nu(\Gamma)&=\int_{\Omega\times\mathbb R^d} \nu^{\mathbf x}(\Gamma)\,d\nu(\mathbf x)
  \\&=\int_{\mathcal{C}} \nu^{\mathbf x}(\Gamma)\,d\nu(\mathbf x)+\int_{\Omega\times \mathbb{R}^d\setminus\mathcal{C}} \nu^{\mathbf x}(\Gamma)\,d\nu(\mathbf x)
\\ &=\int_{\mathcal{C}} \nu^{\mathbf x}(\Gamma\cap \mathcal{C})\,d\nu(\mathbf x)+\int_{\Omega\times \mathbb{R}^d\setminus\mathcal{C}} \nu^{\mathbf x}(\Gamma)\,d\nu(\mathbf x)
\\&\leq \nu(\Gamma\cap\mathcal{C})+\int_{\Omega\times \mathbb{R}^d\setminus\mathcal{C}} \nu^{\mathbf x}(\Gamma)\,d\nu(\mathbf x)
\\&=\int_{\Omega\times \mathbb{R}^d\setminus\mathcal{C}} \nu^{\mathbf x}(\Gamma)\,d\nu(\mathbf x).
 \end{align*}
Therefore, there exists an ergodic measure $\bar{\nu}$ of $\varphi$ such that  $\bar{\nu}(\Gamma)>0$ , and $h_{\bar{\nu}}(\varphi)\ge h_*.$

Without loss of generality,  we  assume that $h_*>0$. Fix $m\ge1$ and set
$$\delta_m:=\min\{m^{-1},h_*/8\}.$$
  Applying \Cref{26-2-10-2156} to the ergodic measure $\bar{\nu}$,  we can choose a finite
 $(\pi_{\mathbb R^d})_*\bar{\nu}$-null boundary measurable partition $\alpha_m$ of $\mathbb R^d$ such that
\begin{equation}
 \label{eq:partition-entropy-lower-bound}
  h_{\bar{\nu}}(\varphi, \alpha_m)
       \ge h_{\bar{\nu}}(\varphi)-\delta_m
       \ge h_*-\delta_m.
\end{equation}
 Now, we  choose
$0<c_m<1/2$  small enough such that
\begin{equation}
 \label{eq:hamming-error}
-c_m\log c_m-(1-c_m)\log(1-c_m)+c_m\log(|\alpha_m|-1)<\delta_m.
\end{equation} 
For $r>0$,  denote
\[
  \mathscr{U}^m_r:=\{x\in \mathbb{R}^d:  B(x,r)\nsubseteq \alpha_m(x)\}.
\]
Since $\bigcap_{r>0}\mathscr{U}^m_r=\partial\alpha_m$ and $\alpha_m$ is a $(\pi_{\mathbb R^d})_*\bar{\nu}$-null boundary measurable partition,  we can choose $r_m>0$ such that
\begin{equation}
 \label{eq:small-boundary-neighborhood}
  (\pi_{\mathbb R^d})_*\bar{\nu}\bigl(\mathscr{U}^m_{3r_m}\bigr)<c_m.
\end{equation}
According to the relative Shannon--McMillan--Breiman theorem (see \cite[Theorem B.0.1]{MR2809170} or \cite[Theorem 3.2]{MET}), we know that for $\bar{\nu}$-a.s. $(\omega,x)\in\Omega\times\mathbb{R}^d$
\begin{equation}
 \label{eq:relative-SMB-component}
\lim_{n\to\infty} -\frac{1}{n}\log\bar{\nu}_\omega\bigl((\alpha_{m})_{\omega}^n(x)\bigr)= h_{\bar{\nu}}(\varphi,\alpha_m).
\end{equation}
 On the other hand,  by Birkhoff's ergodic theorem and \eqref{eq:small-boundary-neighborhood},  for $\bar{\nu}$-a.s. $(\omega,x) \in\Omega\times\mathbb{R}^d$ we have that 
\begin{equation}
 \label{eq:boundary-frequency}
 \lim_{n\to\infty}\frac1n\#\left\{0\le i<n: \varphi(i,\omega)x\in\mathscr{U}^m_{3r_m}\right\}<c_m.
\end{equation}
Let $R_m$ be the full $\bar{\nu}$-measure set on which
\eqref{eq:relative-SMB-component} and
\eqref{eq:boundary-frequency} hold.   It follows that there is a set $\Omega_1$ of full $\mathbb P$-measure such
that  $\bar{\nu}_\omega((R_m)_\omega)=1$ whenever $\omega\in\Omega_1$ and $m\geq 1$.  Define
\begin{equation}
 \label{eq:def-Omega0}
  \Omega_0:=\Omega_1\cap
       \{\omega\in\Omega:
             \bar{\nu}_\omega(\Gamma_\omega)>0\}\subset \pi_{\Omega}(\Gamma).
\end{equation}
Since $\bar{\nu}(\Gamma)>0$, it follows  $\mathbb P(\Omega_0)>0$.   For $N\ge1$, let $R^N_{m}$ be the set of all $(\omega,x)\in R_m$  such that for every $n\ge N$,
\begin{align}
 \label{eq:uniform-SMB-on-RmN}
&\bar{\nu}_\omega\bigl((\alpha_{m})_\omega^n(x)\bigr) \leq \exp\bigl(-n(h_{\bar{\nu}}(\varphi,\alpha_m)-\delta_m)\bigr)
    \le \exp(-n(h_*-2\delta_m)),\\
 \label{eq:uniform-boundary-on-RmN}
 &\sum_{i=0}^{n-1}
   \mathbf{1}_{ \mathscr{U}^m_{3r_m}}
       \bigl(\varphi(i,\omega)x\bigr)
   \le c_m n.
\end{align}
It is obvious that  $R^N_{m}\subset R^{N+1}_{m}$ and $\bigcup_{N\ge1}R^N_{m}=R_m$.  

We claim that there is an integer $N_m^0$ such that  for every $N\ge N_m^0$,
every $(\omega,x)\in R^N_{m}$, and every $n\ge N$,
\begin{equation}
 \label{eq:Bowen-ball-mass-estimate}
\bar{\nu}_\omega\bigl(B(\omega,x,n, 2r_m)\cap(R_{m}^N)_\omega\bigr)\leq \exp(-n(h_*-4\delta_m)).
\end{equation}
Indeed,  assuming that   $\varphi(i,\omega)x\notin\mathscr{U}^m_{3r_m}$ and $y\in B(\omega, x,n, 2r_m)$ for some $0\leq i\leq n-1$,   it is clear that  the points $\varphi(i,\omega)x$ and $\varphi(i,\omega)y$ belong to the same atom of
$\alpha_m$.  By \eqref{eq:uniform-boundary-on-RmN},  we know that 
$$B(\omega,x,n, 2r_m)\subset  \mathbf{B}_{\rho_{(\alpha_m)^n_{\omega}}}(x,c_m).$$
Due to \Cref{26-2-10-2333},  one has that 
$$\#\{ A\in(\alpha_m)^n_{\omega} :  A\cap B(\omega,x,n,2r_m)\neq\varnothing\}\leq\mathrm{B}(c_m,|\alpha_m|, n).$$
 The standard binomial estimate and \eqref{eq:hamming-error} imply that for all
sufficiently large $n$ one has that
\begin{equation}
 \label{eq:number-of-names}
  \mathrm{B}(c_m,|\alpha_m|, n) \leq \exp(2\delta_m n).
\end{equation}
Therefore, we have 
\begin{align*}
\bar{\nu}_\omega\bigl(B(\omega,x,n, 2r_m)\cap(R_{m}^N)_\omega\bigr)&\leq \sum_{\substack{A\in(\alpha_m)^n_{\omega} 
\\ A\cap B(\omega,x,n,2r_m) \cap(R_{m}^N)_\omega\neq\varnothing}}\bar{\nu}_\omega\bigl(A\bigr)
\\\tag{By \eqref{eq:uniform-SMB-on-RmN} and  \eqref{eq:number-of-names} }&\leq \exp(2\delta_m n)\exp(-n(h_*-2\delta_m))
\\&\leq\exp(-n(h_*-4\delta_m)).
\end{align*}

Fix $\omega\in\Omega_0$. For any $m\geq1$  there exists an integer $N=N(\omega,m)\ge N_m^0$ such that $\bar{\nu}_\omega\bigl(
       \Gamma_\omega\cap(R_{m}^N)_\omega\bigr)>0.$  For every sufficiently large $n$, we choose a maximal
$(\omega,n,r_m)$-separated set $E_n(\omega)\subset
       \Gamma_\omega\cap(R_{m}^N)_\omega.$ Then one has that 
\[
 \Gamma_\omega\cap(R_{m}^N)_\omega \subset \bigcup_{x\in E_n(\omega)}B(\omega,x,n, 2r_m).
\]
By using \eqref{eq:Bowen-ball-mass-estimate}, we obtain
\[
\begin{split}
\bar{\nu}_\omega\bigl(
       \Gamma_\omega\cap(R_{m}^N)_\omega\bigr)
 &\le \sum_{x\in E_n(\omega)}
  \bar{\nu}_\omega\bigl(
       B(\omega, x,n, 2r_m)\cap(R_{m}^N)_\omega\bigr)\\
 &\le \#E_n(\omega)\exp\big(-n(h_*-4\delta_m)\big).
\end{split}
\]
This  implies that 
\[
 \frac1n\log s(\omega,n,r_m;\Gamma_\omega)
 \ge h_*-4\delta_m+\frac{1}{n}\log \bar{\nu}_\omega\bigl(
       \Gamma_\omega\cap(R_{m}^N)_\omega\bigr).
\]
Therefore,  one has that 
\[
 \lim_{\epsilon\to0}\liminf_{n\to\infty}
   \frac{\log s(\omega,n,\epsilon;\Gamma_\omega)}{n}
 \ge h_*.
\]

Now, we  give the  proof of  the general case of \Cref{26-5-5-1704}. Define
\[
 \Gamma^+ :=\Big\{\mathbf x\in\Gamma: h_{\nu^{\mathbf{x}}}(\varphi)\ge \nu(\Gamma)^{-1}\int_{\Gamma}h_{\nu^{\mathbf{x}}}(\varphi)\,d\nu(\mathbf{x})\Big\}.
\]
Since $\nu( \Gamma^+)>0,$  we can apply the above argument to $\Gamma^+$ with $h_*:=\nu(\Gamma)^{-1}\int_{\Gamma}h_{\nu^{\mathbf{x}}}(\varphi)\,d\nu(\mathbf{x})$.  Then there exists a positive
$\mathbb P$-measure subset $\Omega_0\subset \pi_{\Omega}(\Gamma^+)$  such that for every $\omega\in\Omega_0$,
\[
 \lim_{\epsilon\to0}\liminf_{n\to\infty}
 \frac{\log s(\omega,n,\epsilon;\Gamma^+_\omega)}{n}\geq\nu(\Gamma)^{-1} \int_{\Gamma}h_{\nu^{\mathbf{x}}}(\varphi)\,d\nu(\mathbf{x}).
\]
Since $\Gamma^+_\omega\subset\Gamma_\omega$ for each $\omega\in\Omega_0$, we can immediately  obtain \Cref{26-5-5-1704}.
\end{proof}

\section{Proof of \Cref{26-7-31-1031}}
\label{26-8-29-2229}
In this section, we combine Pesin theory, Katok's formula, and the general multiplicative mean ergodic theorem in
\Cref{prop:ff-exponential-bound} to prove \Cref{26-7-31-1031}. The main idea is inspired by \cite{MR4562648,MR4081544,MR1619565}.

Throughout this section, we fix an invariant measure
$\nu\in\mathcal I(\varphi)$ and a potential function $f$.  Let  $\tau$ be  a positive integer   such that $f$ is
$\mathscr F_0^{\tau}\otimes\mathscr B(\mathbb R^d)$-measurable, and use the corresponding function $\Lambda_f: (-s_0, s_0)\to\mathbb{R}$  to represent the asymptotic log-moment generating function in \Cref{prop:ff-exponential-bound}.  Let
\[
\nu=\int_{\Omega\times\mathbb R^d}\nu^{\mathbf x}\,d\nu(\mathbf x)
\]
be the ergodic decomposition of $\nu$ for $\Phi_{\tau}$. Since the Wiener measure $\mathbb{P}$ is  $\theta_{\tau}$-ergodic, these components
are also the ergodic measures of the RDS $\varphi_{\tau}$\footnote{Here $\varphi_{\tau}$ denotes the discrete-time RDS over the  $(\Omega,\mathscr{F},\mathbb{P},\theta_{\tau})$, defined by $\varphi_{\tau}(n,\omega,x):=\varphi(n\tau,\omega,x)$ for $n\in\mathbb Z_{+}$. Its associated skew-product mapping  is $\Phi_{\tau}$.}.
Recall that $\mu$ is the Markovian measure of equation
\eqref{26-6-8-1611}. As in the preceding section, write
\[
S_nf:=\sum_{i=0}^{n-1}f\circ\Phi_{i\tau},
\qquad
\|f\|:=\operatorname*{ess\,sup}_{\omega\in\Omega}
\bigl(\|f_\omega\|_\infty+\operatorname{Lip}(f_\omega)\bigr).
\]

\subsection{Uniform estimates of Pesin blocks}
\label{26-9-5-1342}

The purpose of this subsection is to establish local estimates of neutralized Bowen balls. Recall that $\Gamma$ is the Lyapunov regular
set in \Cref{26-6-13-1113}. For $1\leq k\leq d$, set
\begin{equation}\label{eq:Zk}
\Gamma_k:=\{\mathbf x\in\Gamma:
\chi_\nu^k(\mathbf x)>0\geq\chi_\nu^{k+1}(\mathbf x)\},
\end{equation}
where $\chi_\nu^{d+1}(\mathbf x):=0$, and also set
\[
\Gamma_0:=\{\mathbf x\in\Gamma:\chi_\nu^1(\mathbf x)\leq0\}.
\]
Thus $\Gamma=\bigcup_{k=0}^d\Gamma_k$ modulo $\nu$, and each
$\Gamma_k$ is $\Phi_{\tau}$-invariant modulo $\nu$. Let $\operatorname{vol}$ denote Lebesgue measure on $\mathbb R^d$.
For $\epsilon>0$ and $n\geq1$, write
\begin{align*}
B_{\tau}(\omega,x,n,\epsilon)&:=B(\omega,x,n\tau,\epsilon)
\\&=\left\{y\in\mathbb R^d:|\varphi(i,\omega,y)-\varphi(i,\omega,x)|\leq \epsilon \text{ for }i=0,\ldots,n\tau-1\right\}.
\end{align*}

\begin{lemma}
\label{lem:pointwise-shrinking-bowen}
For $\nu$-almost every $\mathbf x=(\omega,x)\in\Gamma$, one has
\begin{equation}\label{26-8-23-1524}
\lim_{\epsilon\to0}\limsup_{n\to\infty}
-\frac1n\log\operatorname{vol}B_{\tau}(\omega,x,n,e^{-n\epsilon})\leq \tau\chi_\nu^+(\mathbf x),
\end{equation}
where $\chi_\nu^+(\mathbf x)=\sum_{i=1}^d\chi_\nu^i(\mathbf x)^+$.
Moreover, for every $\epsilon>0$, we have 
\begin{equation}\label{26-10-3-1426}
\lim_{n\to\infty}
\sup_{y\in B_{\tau}(\omega,x,n,e^{-n\epsilon})}\left|\frac1nS_nf(\omega,y)-\mathbb E_{\nu^{\mathbf x}}[f]\right|=0.
\end{equation}
\end{lemma}

\begin{proof}
Recall that $\nu(\mathfrak A)=1$. For
$\mathbf x\in\Gamma_k$, let
\[
E^u(\mathbf x):=\bigoplus_{\lambda_\nu^i(\mathbf x)>0}E^i(\mathbf x),
\qquad
E^{cs}(\mathbf x):=\bigoplus_{\lambda_\nu^i(\mathbf x)\leq0}E^i(\mathbf x).
\]
When $k\geq1$,  we choose a unit decomposable $k$-vector $v^u(\mathbf x)$ spanning $E^u(\mathbf x)$. The multiplicative
ergodic theorem implies that
\begin{equation}\label{eq:MET-unstable-jacobian}
\lim_{n\to\infty}\frac1n
\log\left|
\bigl(\mathrm d_x\varphi(n\tau,\omega)\bigr)^{\wedge k}
 v^u(\mathbf x)\right|=\tau\chi_\nu^+(\mathbf x).
\end{equation}
Together with the integrability estimates in \Cref{26-6-9-1713},
the neutralized Bowen-volume argument of \cite[Proposition~5.11]{Lianzhang} yields \eqref{26-8-23-1524}. When $k=0$, \eqref{26-8-23-1524} obviously holds by  \cite[Proposition~5.11]{Lianzhang}.

Since $f$ is a potential function, for $\mathbb{P}$-a.s. $\omega\in\Omega$ and any $x\in\mathbb{R}^d$, we have
\begin{equation}\label{eq:bowen-potential-distortion}
\begin{aligned}
|S_nf(\omega,y)-S_nf(\omega,x)|&\leq\|f\|\sum_{i=0}^{n-1}|\varphi(i\tau,\omega,y)-\varphi(i\tau,\omega,x)|\leq n\|f\|e^{-n\epsilon},
\end{aligned}
\end{equation}
whenever $y\in B_{\tau}(\omega,x,n,e^{-n\epsilon})$. By Birkhoff's ergodic theorem for $\Phi_{\tau}$, we have for $\nu$-a.s. $\mathbf x=(\omega,x)\in\Omega\times\mathbb{R}^d$,
\[
\lim_{n\to\infty}\frac1nS_nf(\omega,x)
=\mathbb E_{\nu^{\mathbf x}}[f].
\]
Combining this with \eqref{eq:bowen-potential-distortion}, we obtain \eqref{26-10-3-1426}.
\end{proof}

We now construct the Pesin blocks used in the next subsection.

\begin{proposition}[Pesin blocks]
\label{lem:nonergodic-pesin-blocks}
Fix $0\leq k\leq d$ and $\delta>0$. Then
\[
\Gamma_k=
\bigcup_{\mathbf I\in\mathbb Q_{+}^2\times\mathbb Q\times\mathbb N^3}
\Gamma_k(\mathbf I;\delta)\pmod\nu,
\]
where $\mathbf I=(h,\chi,\mathbb A,R,q,N)$, and
$\Gamma_k(\mathbf I;\delta)$ consists of the points
$\mathbf x=(\omega,x)\in\Gamma_k$ satisfying $|x|\leq R$, and 
\begin{equation}\label{eq:block-bins}
\begin{gathered}
h\leq h_{\nu^{\mathbf x}}(\varphi_{\tau})<h+\delta,
\quad \chi\leq \chi_\nu^+(\mathbf x)<\chi+\delta,\quad\mathbb A\leq\mathbb E_{\nu^{\mathbf x}}[f]<\mathbb A+\delta,
\end{gathered}
\end{equation}
and, for every $n\geq N$,
\begin{align}
&\operatorname{vol}B_{\tau}(\omega,x,n,e^{-n/q})
\geq\exp\bigl(-n\tau(\chi+2\delta)\bigr),
\label{eq:block-core-volume}\\
&\sup_{y\in B_{\tau}(\omega,x,n,e^{-n/q})}
\left|\frac1nS_nf(\omega,y)-\mathbb A\right|\leq2\delta.
\label{eq:block-core-potential}
\end{align}
\end{proposition}

\begin{proof}
The rational bins in \eqref{eq:block-bins} cover all finite values of the time-$\tau$ entropy, the  positive-exponent sum, and the
Birkhoff average. By \Cref{lem:pointwise-shrinking-bowen}, for
almost every point one may choose an integer $q$ sufficiently large and then $N$ sufficiently large so that \eqref{eq:block-core-volume}
holds for every $n\geq N$. Increasing $N$ also gives \eqref{eq:block-core-potential}, by \eqref{26-10-3-1426} and the bin for $\mathbb E_{\nu^{\mathbf x}}[f]$. This proves \Cref{lem:nonergodic-pesin-blocks}.
\end{proof}

\subsection{Proof of \Cref{26-7-31-1031}}

In this subsection, we combine Katok's formula in \Cref{26-5-5-1704} with the uniform Pesin-block estimates above to relate asymptotic log-moment generating function, metric entropy, and Lyapunov exponents. Then we complete the proof of \Cref{26-7-31-1031} if the invariant measure satisfies  entropy formula. 

\begin{proposition}
\label{26-8-14-1929}
Let $f$ be a potential function and $\nu\in\mathcal I(\varphi)$. For every $0\leq k\leq d$ with $\nu(\Gamma_k)>0$ and  any $|s|<s_0$, we have
\begin{equation}\label{eq:finite-block-average}
\Lambda_f(s)\geq\frac1{\nu(\Gamma_k)}
\int_{\Gamma_k}\left(h_{\nu^{\mathbf x}}(\varphi_{\tau})-\tau\chi_\nu^+(\mathbf x)+s\int_{\Omega\times\mathbb{R}^d} f\,d\nu^{\mathbf x}\right)\,d\nu(\mathbf x).
\end{equation}
\end{proposition}

\begin{proof}
It is enough to prove the stronger estimate
\begin{equation}\label{eq:vlb-essential-supremum}
\Lambda_f(s)\geq\operatorname*{ess\,sup}_{\mathbf x\in\Gamma_k}\Big(h_{\nu^{\mathbf x}}(\varphi_{\tau})-\tau\chi_\nu^+(\mathbf x)
+s\int_{\Omega\times\mathbb{R}^d} f\,d\nu^{\mathbf x}\Big),
\end{equation}
where the essential supremum is taken with respect to $\nu$.
Fix $\delta>0$, a constant  $s\in(-s_0, s_0)$, and a number $a$ strictly
smaller than the essential supremum in
\eqref{eq:vlb-essential-supremum}. Then
\[
\mathcal{S}=\Big\{\mathbf x\in\Gamma_k: h_{\nu^{\mathbf x}}(\varphi_{\tau})-\tau\chi_\nu^+(\mathbf x)
+s\int_{\Omega\times\mathbb{R}^d} f\,d\nu^{\mathbf x}>a\Big\}
\]
has positive $\nu$-measure. By \Cref{lem:nonergodic-pesin-blocks},
there is
\[
\mathbf I_0=(h_0,\chi_0,\mathbb A_0,R_0,q_0,N_0)
\in\mathbb Q_{+}^2\times\mathbb Q\times\mathbb N^3
\]
such that $\nu(\mathcal{S}\cap\Gamma_k(\mathbf I_0;\delta))>0$. By regularity of $\nu$, we can choose a compact subset $G$ of  $\mathcal{S}\cap\Gamma_k(\mathbf I_0;\delta)$ with  positive $\nu$-measure. In particular, for every $\mathbf x\in G$,
\begin{equation}\label{eq:vlb-block-bins}
h_0\leq h_{\nu^{\mathbf x}}(\varphi_{\tau})<h_0+\delta,
\end{equation}
and the other estimates  in \eqref{eq:block-bins} imply that 
\begin{equation}\label{eq:vlb-bin-comparison}
h_0-\tau\chi_0+s\mathbb A_0>a-(1+|s|+\tau)\delta.
\end{equation}

Apply \Cref{26-5-5-1704} to $G$ for the RDS $\varphi_{\tau}$,  there is a measurable set
$\widehat\Omega\subseteq\pi_\Omega(G)$ of positive
$\mathbb P$-measure such that, for every $\omega\in\widehat\Omega$ one has 
\begin{equation}\label{eq:vlb-local-entropy}
\lim_{\epsilon\to 0}\liminf_{n\to\infty}
\frac1n\log s(\omega,n\tau,\epsilon;G_\omega)\geq\frac1{\nu(G)}\int_Gh_{\nu^{\mathbf x}}(\varphi_{\tau})\,d\nu(\mathbf x)
\geq h_0,
\end{equation}
where we use the fact that  every separated set of length $n$ for $\varphi_{\tau}$ is also separated at length $n\tau$ for the original RDS $\varphi$.
For integers $q,N\geq1$, define
\[
\Omega(q,N):=\left\{\omega\in\widehat\Omega: s(\omega,n\tau,1/q;G_\omega)\geq e^{n(h_0-\delta)}
\text{ for every }n\geq N\right\}.
\]
These sets are $\mathscr{F}$-measurable by  compactness of  $G$  and  Remark \ref{26-9-5-1802}.  By \eqref{eq:vlb-local-entropy}, we have 
\[
\widehat\Omega=\bigcup_{q,N\geq1}\Omega(q,N).
\]
Thus there exists some $\Omega_1:=\Omega(q_1,N_1)$ with positive $\mathbb P$-measure. For any $\omega\in\Omega_1$ and
$n\geq\max\{N_0,N_1\}$, we can choose an $(\omega,n\tau,1/q_1)$-separated set $E_n(\omega)\subseteq G_\omega$
such that
\begin{equation}\label{eq:vlb-separated-cardinality}
\#E_n(\omega)\geq e^{n(h_0-\delta)}.
\end{equation}
Let $U:=\{x\in\mathbb R^d:|x|\leq R_0+1\}$. For each $\omega\in\Omega_1$ and  all sufficiently large $n$, we have
\begin{align*}
\int_Ue^{sS_nf(\omega,y)}\,d\operatorname{vol}(y)&\geq\sum_{x\in E_n(\omega)}\int_{B_{\tau}\big(\omega,x,n, \frac{1}{2q_1}\big)}
e^{sS_nf(\omega,y)}\,d\operatorname{vol}(y)\\
&\geq\sum_{x\in E_n(\omega)}
\int_{B_{\tau}(\omega,x,n,e^{-n/q_0})}
e^{sS_nf(\omega,y)}\,d\operatorname{vol}(y)\\
&\geq\#E_n(\omega)\,
e^{-n\tau(\chi_0+2\delta)}e^{n(s\mathbb A_0-2|s|\delta)}\\
&\geq e^{n(h_0-\tau\chi_0+s\mathbb A_0-(1+2\tau+2|s|)\delta)}\\
&\geq e^{n(a-(2+3|s|+3\tau)\delta)}.
\end{align*}
Integrating the above inequality over $\Omega_1$, we can obtain 
\begin{equation}\label{eq:vlb-annealed-lower}
\liminf_{n\to\infty}\frac1n\log
\mathbb E\int_Ue^{sS_nf(\omega,y)}\,d\operatorname{vol}(y)
\geq a-(2+3|s|+3\tau)\delta.
\end{equation}
By the locally uniform convergence in  \Cref{prop:ff-exponential-bound} and Fubini's  theorem, one has that
\begin{equation}\label{eq:vlb-identification}
\lim_{n\to\infty}\frac1n\log
\mathbb E\int_Ue^{sS_nf(\omega,y)}\,d\operatorname{vol}(y)=\Lambda_f(s).
\end{equation}
Combining \eqref{eq:vlb-annealed-lower} and
\eqref{eq:vlb-identification}, we obtain \eqref{eq:vlb-essential-supremum} by letting $\delta$ decrease to zero,  and then letting $a$ increase to the essential supremum.
\end{proof}

We now use the entropy formula and the derivative of $\Lambda_f$ at zero to complete the proof of \Cref{26-7-31-1031}.

\begin{proof}[Proof of \Cref{26-7-31-1031}]
Let $\nu\in\mathcal I(\varphi)$  satisfy the entropy formula \eqref{26-7-31-1028}. Multiplying \eqref{eq:finite-block-average} by $\nu(\Gamma_k)$ and summing over the sets of positive measure, including $k=0$,  we have 
\begin{align*}
\Lambda_f(s)&\geq\int_\Gamma\left(
h_{\nu^{\mathbf x}}(\varphi_{\tau})-\tau\chi_\nu^+(\mathbf x)+s\int_{\Omega\times\mathbb{R}^d} f\,d\nu^{\mathbf x}\right)\,d\nu(\mathbf x)\\
&=h_\nu(\varphi_{\tau})-\tau\chi_\nu^++s\mathbb E_\nu[f]\\
&=\tau\bigl(h_\nu(\varphi)-\chi_\nu^+\bigr)+s\mathbb E_\nu[f]\\
&=s\mathbb E_\nu[f].
\end{align*}
Here we used the fact $h_\nu(\varphi_{\tau})=\tau h_\nu(\varphi)$, and the entropy formula.
Since $\Lambda_f(0)=0$, dividing by $s>0$ and $s<0$, respectively,
and letting $s\to0$ yields
\begin{equation}\label{26-10-3-1605}
\mathbb E_\nu[f]=\Lambda_f'(0)=\mathbb E_\mu[f].
\end{equation}
Since $f$ was arbitrary, this identity holds for every potential
function.

\begin{claim}\label{26-10-3-1622}
The measures $\mu$ and $\nu$ coincide on the future. Namely,
$\mu(A)=\nu(A)$ for every
$A\in\mathscr F_0^\infty\otimes\mathscr B(\mathbb R^d)$.
\end{claim}

\begin{proof}
Fix an integer $q\geq1$ and $C\in\mathscr F_0^q$. Define finite Borel measures on $\mathbb R^d$ by
\[
\nu_C(B):=\nu(C\times B),
\qquad\mu_C(B):=\mu(C\times B),
\]
for  $B\in\mathscr B(\mathbb R^d).$ For every $g\in\operatorname{Lip}_b(\mathbb R^d)$,
$f(\omega,x)=\mathbf1_C(\omega)g(x)$ is a potential function. Thus by \eqref{26-10-3-1605} we have 
\[
\int_{\mathbb R^d}g\,d\nu_C
=\int_{\mathbb R^d}g\,d\mu_C.
\]
Then we immediately  obtain  $\nu_C=\mu_C$ by  using  the dual Lipschitz distance
\cite[Chapter~1.2.3]{MR3443633}.

Since the $\sigma$-algebras $\mathscr F_0^q$ are increasing,
\[
\mathscr A:=\bigcup_{q\geq1}\mathscr F_0^q
\]
is an algebra, with $\sigma(\mathscr A)=\mathscr F_0^\infty$.
Hence
\[
\mathscr P:=\{C\times B:C\in\mathscr A,\ B\in\mathscr B(\mathbb R^d)\}
\]
is a $\pi$-system containing $\Omega\times\mathbb R^d$ and
generating $\mathscr F_0^\infty\otimes\mathscr B(\mathbb R^d)$.
Since $\nu$ and $\mu$ agree on $\mathscr P$, the $\pi$-$\lambda$ theorem implies \Cref{26-10-3-1622}.
\end{proof}

Finally, fix a constant  $\ell\geq1$, a set $C\in\mathscr F_{-\ell}^\infty$, and a bounded Borel function
$g:\mathbb R^d\to\mathbb R$. It is clear that the function $\mathbf1_C\circ\theta_\ell$ is
$\mathscr F_0^\infty$-measurable, and $(\omega,x)\mapsto\varphi(\ell,\omega,x)$ is
$\mathscr F_0^\ell\otimes\mathscr B(\mathbb R^d)$-measurable.
By invariance and \Cref{26-10-3-1622},
\begin{align*}
\int_{\Omega\times\mathbb R^d}\mathbf1_C(\omega)g(x)\,d\nu(\omega,x)
&=\int_{\Omega\times\mathbb R^d}\mathbf1_C(\theta_\ell\omega)
g(\varphi(\ell,\omega,x))\,d\nu(\omega,x)\\
&=\int_{\Omega\times\mathbb R^d}\mathbf1_C(\theta_\ell\omega)
g(\varphi(\ell,\omega,x))\,d\mu(\omega,x)\\
&=\int_{\Omega\times\mathbb R^d}\mathbf1_C(\omega)g(x)\,d\mu(\omega,x).
\end{align*}
Since $\sigma(\bigcup_{\ell\geq1}\mathscr F_{-\ell}^\infty)=\mathscr F$,
another monotone class argument yields $\nu=\mu$.
Thus every invariant measure satisfying the entropy formula coincides with $\mu$. This completes the proof of \Cref{26-7-31-1031}.
\end{proof}

\section{SRB bifurcation for the Hopf normal form}
\label{26-8-29-2233}

For  2D additive-noise Hopf normal form,  it  has been  investigated a lot for its sample top Lyapunov exponents  \cite{DoanEtAl2018,ChemnitzEngel2023,Baxendale2025}. In this section, we mainly verify that the assumptions  \hyperref[26-6-10-1501]{\textnormal{(H1)}}--\hyperref[26-6-10-1503]{\textnormal{(H3)}} hold.  Together with results in  \cite{breden2024rigorous} and the continuity of top Lyapunov exponents with respect to the shear parameter, we establish the bifurcation of the existence of SRB measures for   2D additive-noise stochastic Hopf normal form.

For $a,\sigma>0$ and $\alpha,\beta,b\in\R$, the Hopf normal form is considered as following 
\begin{equation}
\label{eq:hopf-vector-form}
 dX_t=\Big((\alpha-a|X_t|^2)X_t+(\beta-b|X_t|^2)J X_t\Big)\, dt+\sigma e_1\,dW^1_t+\sigma e_2\,dW^2_t,
\end{equation}
where   $e_1=(1,0)^{T}$, $e_2=(0,1)^{T}$ and 
\[
 X_t=\begin{pmatrix}
 X^{(1)}_t
\\ X^{(2)}_t
\end{pmatrix},
 \qquad
 J=\begin{pmatrix}0&-1\\1&0\end{pmatrix}.
\]
Here $a$ controls radial contraction and $b$ is the shear coefficient. Since $\langle x ,Jx\rangle=0$, there exists $C_{\alpha,a}>0$ such that
\begin{equation}
\label{eq:hopf-radial-dissipation}
 \bigl\langle x,F_{}(x)\bigr\rangle=\alpha|x|^2-a|x|^4\leq-\frac a2|x|^4+C_{\alpha,a}.
\end{equation}
where $F(x)= (\alpha-a|x|^2)x+(\beta-b|x|^2)Jx.$ The  first two derivatives of drift term   have
polynomial growth, and the diffusion matrix is $\sigma^2\text{Id}$.  Thus \hyperref[26-6-10-1501]{\textnormal{(H1)}}--\hyperref[26-6-10-1503]{\textnormal{(H3)}} hold with $p=4$.  Therefore, it generates a $C^2$ RDS and has a unique Markovian measure  $\mu_{\alpha,\beta,a,b,\sigma}$.  We write its corresponding  top Lyapunov exponent as $\chi^1(\alpha,\beta,a,b,\sigma)$.

The following theorem describes the transition from nonexistence to existence of SRB measures as the shear parameter varies.
A key feature is that it  resolves the critical case where the top sample Lyapunov exponent vanishes. To the best of our knowledge, the nonexistence of SRB measures for the stochastic Hopf normal form (even for other stochastic models) in this critical case has not previously been established.

\begin{thm}
\label{cor:hopf-srb-transition-continuity}
Let $a,\sigma>0$, $\beta\in\R$ and $\alpha=\sigma\sqrt{2a}$. There exist constants  $s_-<c_-\le c_+<s_+$ such that 
\begin{equation}
\label{eq:hopf-original-strict-dichotomy}
 \begin{cases}
 \chi^1(\alpha,\beta,a,b,\sigma)<0,
  &\text{ when }b\in [s_-, c_-),\\[1mm]
  \chi^1(\alpha,\beta,a,b,\sigma)=0,
  &\text{ when }b\in\{c_-, c_+\},\\[1mm]
 \chi^1(\alpha,\beta,a,b,\sigma)>0,
  &\text{ when }b\in(c_+, s_+].
 \end{cases}
 \end{equation}
 Particularly, equation \eqref{eq:hopf-vector-form} admits no SRB measure for $b\in[s_-,c_-]\cup\{c_+\}$. For any $b\in(c_+,s_+]$, its Markovian measure $\mu_{\alpha,\beta,a, b,\sigma}$ is an SRB measure.
\end{thm}

\begin{proof}
 We first  recall the following symmetry and scaling properties of the top Lyapunov exponent of equation \eqref{eq:hopf-vector-form}. The detail proof can refer to \cite[Section~1.1 and Proposition~3.1]{Baxendale2025}. For any $a,\sigma>0$ and $\alpha,\beta,b\in\R$, the top Lyapunov exponent  satisfies
\begin{equation}
\label{eq:hopf-direct-reference-scaling}
 \chi^1(\alpha,\beta,a,b,\sigma)
 =\frac{\sigma\sqrt a}{2\sqrt2}\,
\chi^1\left(
   \frac{2\sqrt2\,\alpha}{\sigma\sqrt a},
   0,4,\frac{4|b|}{a},\sqrt2
 \right).
\end{equation}
According to \cite[Theorem~4 and Corollary~5]{breden2024rigorous}, one has
\begin{equation}
\label{eq:hopf-reference-endpoint-signs}
 \chi^1(4,0,4,21.5322,\sqrt2)<0
 \quad\text{and}\quad
 \chi^1(4,0,4,21.5381,\sqrt2)>0.
 \end{equation}
   Let $\alpha=\sigma\sqrt{2a}$. By using \eqref{eq:hopf-direct-reference-scaling} and  \eqref{eq:hopf-reference-endpoint-signs}, one has that
\begin{equation}
\label{eq:hopf-original-endpoint-signs}
\chi^1(\alpha,\beta,a, s_-,\sigma)<0
 <\chi^1(\alpha,\beta,a, s_+,\sigma),
\end{equation}
where $s_-:=5.38305\,a$ and $s_+:=5.384525\,a$. Next, we claim that the top Lyapunov exponent is continuous with respect to the shear parameter. 
\begin{claim}
\label{lem:hopf-shear-continuity}
For fixed $\alpha,\beta\in\mathbb R$ and $a,\sigma>0$, the map
$b\mapsto \chi^1(\alpha,\beta,a,b,\sigma)$ is continuous.
\end{claim}
\begin{proof}
Recall that $\varphi$ is the RDS generated by  equation \eqref{eq:hopf-vector-form}.  We first  prove that the direction process 
$$w_t:=(\varphi(t, \omega,x), [\mathrm{d}_x\varphi(t,\omega)v])$$ 
on  $\mathcal M=\R^2\times\mathbb P^1$ has a unique stationary measure. For convenience, projective directions are represented locally by unit vectors. For any $(x,[v])\in\mathcal{M}$,  $x$ and $v$ is called  horizontal component and vertical component, respectively.  Let $H_e(\cdot)\equiv (e,0)$ be  the corresponding constant horizontal vector field  for $e\in\mathbb R^2$,   and  let $z^{\operatorname{vert}}$ be the vertical component for $z\in\mathcal{M}$.  Usually,   we have following identification
$$T_{[v]}\mathbb P^1\simeq v^\perp=\operatorname{span}\{Jv\}.$$
By a direct calculation,   the direction process satisfies the Stratonovich equation
\[
 dw_t=\widehat F(w_t)\,dt+\sigma H_{e_1}(w_t)\circ dW_t^1+\sigma H_{e_2}(w_t)\circ dW_t^2,
\]
where 
\[ \widehat F(x,[u]):=\left(F(x),\Pi_{u^\perp}\nabla F(x)u\right),
\] and $\Pi_{u^\perp}=\text{Id}-uu^{T}$ is  the orthogonal projection onto $u^\perp$.

Fix  a  point $z_0=(x_0,[v_0])\in\mathcal{M}$,  and  let $\eta_0=Jv_0$.   It is clear that $H_{e_1}$ and $H_{e_2}$ span the entire horizontal subspace.  Since $H_{v_0}$ and $H_{\eta_0}$ are constant horizontal fields,  one has that
\[
[H_{v_0},\widehat{F}]^{\mathrm{vert}}(y,[u])
=
\Pi_{u^\perp}\nabla^2F(y)[v_0,u],
\]
and then
\begin{equation}
\label{eq:hopf-double-bracket-general}
[H_{\eta_0},[H_{v_0},\widehat{F}]]^{\mathrm{vert}}(z_0)=\Pi_{v_0^\perp}\nabla^3F(x_0)[\eta_0,v_0,v_0]=-2a\,\eta_0
\neq0.
\end{equation}
 Thus, the double Lie bracket products a non-zero  vertical direction in $T_{[v_0]}\mathbb P^1$. By  H\"ormander's theorem,  it follows that the direction process $w_t$ is strong Feller.

 For a Cameron--Martin path $h$, the corresponding controlled system is
\[
 \dot x=F(x)+\sigma\dot h,
 \qquad
 \dot v=\Pi_{v^\perp}\nabla F(x)v,
\]
where $|v|=1.$  Every prescribed piecewise smooth base path $x$ can be realized by choosing
\[
 \dot h(t)=\sigma^{-1}\bigl(\dot x(t)-F(x(t))\bigr).\]
For $z:=(x,v)\in\mathcal{M}$, we can write 
 $$x=r e_{\phi}:=r(\cos\phi,\sin\phi) \quad\text{and}\quad v=e_{\theta}:=(\cos\theta,\sin\theta).$$
 Then  direction process satisfies  
\begin{equation}
 \label{eq:hopf-controlled-direction}
 \dot\theta
 =\beta-r^2\Gamma_{a,b}(\phi-\theta),
  \end{equation}
where $\Gamma_{a,b}(\delta):=b+a\sin(2\delta)+2b\cos^2\delta.$ Fix arbitrary initial  and terminal states   $z_0=(x_0,[e_{\theta_0}])$ and $ z_1=(x_1,[e_{\theta_1}])$. We will  construct a controlled trajectory from $z_0$ to $z_1$ in time $1$.  Firstly, we  choose  smooth base paths
$x_-:[0, 1/3]\to\mathbb R^2$ and $x_+:[0, 1/3]\to\mathbb R^2$ satisfying
\[
 x_-(0)=x_0,\qquad x_-(1/3)=\textbf{0},
 \qquad
 x_+(0)=\textbf{0},\qquad x_+(1/3)=x_1.
\]
Let $[e_{\theta_a}]$  be the direction at time $1/3$ along $x_-$.  Equation \eqref{eq:hopf-controlled-direction} along the  path $x_+$ induces
a diffeomorphism
$\Psi:\mathbb P^1\to\mathbb P^1$. Let $\theta_b\in\mathbb R$ such that $\Psi([e_{\theta_b}])=[e_{\theta_1}]$.  Therefore, we only need to  construct a middle base loop of length $1/3$, based at the origin, which carries $[e_{\theta_a}]$ to $[e_{\theta_b}]$.

It is clear that we can choose a $\delta_0$ such that
$\Gamma_0:=\Gamma_{a,b}(\delta_0)\ne0$.
Let $\bar r\in C_c^\infty((0, 1/3))$ be nonnegative and not identically zero, and  we can choose some  $n\in\mathbb Z$ such that
\[
 R^2
 :=\frac{\theta_a+\beta/3-\theta_b-n\pi}
          {\Gamma_0 I_0}>0,
\]
where  $I_0:=\int_0^{1/3}\bar r(t)^2\,dt>0.$  Define
\[
 \theta_R(t)
 :=\theta_a+\beta t
   -R^2\Gamma_0\int_0^t\bar r(s)^2\,ds,
 \qquad
 x_R(t):=R\bar r(t)e_{\theta_R(t)+\delta_0}.
\]
It is not hard to see that $(x_R,\theta_R)$ satisfies the equation \eqref{eq:hopf-controlled-direction}. Let $x_R(0)=x_R(1/3)=\textbf{0}$, then we have 
\[
 \theta_R(0)=\theta_a,
 \qquad
 \theta_R(1/3)=\theta_b+n\pi.
\]
Thus this middle loop has the required initial and terminal projective directions. Concatenating $x_-$, $x_R$, and $x_+$ , we obtain  a piecewise smooth base path on $[0,1]$ whose associated controlled trajectory starts at $z_0$ and ends at $z_1$.
The control defined above is a Cameron--Martin path.
By Stroock--Varadhan support theorem (see \cite[Theorem~19.1]{MR2604669}), after a standard localization around the controlled trajectory, we know that the direction process $w_t$ is irreducible. Hence, this process admits a unique stationary measure.

Let   $F_b$ be  the drift term  in \eqref{eq:hopf-vector-form}, and let $\hat{\gamma}_b$ be  the unique stationary measure of the direction process on $\mathcal M=\mathbb R^2\times\mathbb P^1$. By \cite[Section~1]{Baxendale2025}, we know that the unique stationary measure of its spatial process is 
\[
 \gamma(dx)=Z^{-1}\exp\!\left(\frac{\alpha|x|^2-\frac a2|x|^4}{\sigma^2}\right)\,dx,
\]
where $Z$ is the normalized constant. 

Let $\{b_n\}_{n\in\mathbb{N}}$ be a sequence which converges to $b$. Since $\mathbb P^1$ is compact and all $\hat{\gamma}_{b_n}$ have  marginal $\gamma$, the sequence $\{\hat{\gamma}_{b_n}\}$ is tight. Local continuous dependence of the coefficients on $b$, and the uniform radial estimate \eqref{eq:hopf-radial-dissipation} imply  that the corresponding transition semigroups converge locally uniformly for  bounded continuous test functions.  By passing to the limit in the invariance identity, one has  that every weak limit of $\{\hat{\gamma}_{b_n}\}$ is  a stationary measure for the process with parameter $b$. The reader can refer to the same argument in \cite[Lemmas~3.9 and~3.12]{ChemnitzEngel2023}. Due  to the uniqueness  of stationary measure,  we have  $\lim_{n\to\infty}\hat{\gamma}_{b_n}=\hat{\gamma}_b$.

By the Furstenberg--Khasminskii formula (see \cite[Section~1.1]{Baxendale2025}),
\[
 \chi^1(\alpha,\beta,a,b,\sigma) =\int_{\mathcal M}\mathcal{G}_b\,d\hat{\gamma}_b,
\]
where  $\mathcal{G}_b(x,[v]):=\frac{\langle \nabla F_b(x)v,v\rangle}{|v|^2}$.  Note that  $\lim_{n\to\infty}\mathcal{G}_{b_n}= \mathcal{G}_b$ converges locally uniformly, and
$$|\mathcal{G}_{b_n}(x,[v])|\le C(1+|x|^2)$$
for some constant $C$ independent of  $n$.  Consequently,  one has that 
$$\lim_{n\to\infty}\int_{\mathcal M}\mathcal{G}_{b_n}\,d\hat{\gamma}_{b_n}=\int_{\mathcal M}\mathcal{G}_b\,d\hat{\gamma}_b,$$
which  proves the claim.
\end{proof}

Letting  
\begin{align*}
&c_-:=\min\{b>s_-:  \chi^1(\alpha,\beta,a,b,\sigma)=0\} 
\\& c_+:=\max\{b<s_+:  \chi^1(\alpha,\beta,a,b,\sigma)=0\},
\end{align*}
we  obtain \Cref{cor:hopf-srb-transition-continuity} by \Cref{lem:hopf-shear-continuity} and  \Cref{26-8-25-1504}.
\end{proof}

\subsection*{Acknowledgments and Statements}
The authors would like to thank  Rongchang Liu for valuable comments on an earlier version of this manuscript.  Zeng Lian was supported by  the NSFC of China  (Grant  No. 12090010) and National Tianyuan Mathematics Southwest Center (Grant  No. 12226102).  Jianhua Zhang was supported by  the China Postdoctoral Science Foundation (Grant  No. 2025M774190). This work is original research and has not been submitted to any other journals for publication. We declare that we do not have any commercial or associative interest that represents a conflict of interest in connection with the work submitted.

\subsection*{Declaration of AI assistance}

The mathematical ideas and proof strategy underlying the initial
version of this work were developed by the authors. In that
version, we used moment Lyapunov exponents to establish the
variational identity
\[
 h_{\mu}(\varphi)
 =
 \sup_{\nu\in\mathcal{I}_{\mathrm{ent}}(\varphi)}
 h_{\nu}(\varphi),
\]
where $\mu$ is the Markovian invariant measure and
$\mathcal{I}_{\mathrm{ent}}(\varphi)$ denotes the collection of
invariant measures satisfying the entropy formula.

After completing these proofs, we uploaded the manuscript to Codex, using the GPT-6-Astra model, and asked whether the
uniqueness of SRB measures could be established. Building on our original approach, the AI developed an argument for
uniqueness within the framework of the multiplicative mean ergodic theorem now presented in \Cref{26-8-17-1609}.
This formulation is closely related to the Feynman-Kac approach used for moment Lyapunov exponents in the initial version and follows the same underlying spectral strategy. We subsequently revised and reorganized the manuscript to incorporate this argument and present the
results in their current form.

AI assistance was also used to discuss technical details of the proofs and improve the exposition during revision. The authors
take full responsibility for the content of this paper.
\providecommand{\bysame}{\leavevmode\hbox to3em{\hrulefill}\thinspace}
\providecommand{\MR}{\relax\ifhmode\unskip\space\fi MR }
\providecommand{\MRhref}[2]{%
  \href{http://www.ams.org/mathscinet-getitem?mr=#1}{#2}
}
\providecommand{\href}[2]{#2}

\end{document}